\documentclass[11pt,oneside,english]{amsart}

\usepackage{amsmath,amsthm}
\usepackage{mathtools}
\mathtoolsset{mathic}
\usepackage{empheq}[2004/10/10]
\usepackage{color}
\usepackage{tikz}

\usepackage{microtype}

\usepackage{ifxetex,ifluatex}
\ifxetex
\usepackage{xltxtra}
\else
\ifluatex
\usepackage{luatextra}
\else
\usepackage{amssymb}
\usepackage[mathscr]{eucal}
\usepackage[utf8]{inputenc}
\usepackage[T1]{fontenc}
\usepackage{mathpazo,tgcursor,tgheros,tgpagella}
\usepackage{textcomp}

\csname else\expandafter\endcsname
\romannumeral -`0
\fi
\fi
\iftrue\relax%
\defaultfontfeatures{Ligatures=TeX}
\usepackage[math-style=TeX, vargreek-shape=TeX]{unicode-math}
\AtBeginDocument{%
  \let\setminus\smallsetminus
}
\fi

\usepackage[letterpaper,margin=1in]{geometry}

\usepackage[footnotesize,nooneline,bf]{caption}
\usepackage{enumitem}
\newlist{enumerate*}{enumerate*}{1}
\setlist[enumerate*]{label=(\arabic*),
  itemjoin={{, }}, itemjoin*={{, and }}}
\setenumerate{label=(\arabic*)}
\setitemize{label=\tiny\textbullet}

\usepackage{xspace}
\usepackage{array,hhline}

\usepackage[unicode=true, pdfusetitle,
 bookmarks=true,bookmarksnumbered=false,bookmarksopen=false,
 breaklinks=false,pdfborder={0 0 0},backref=false,colorlinks=false]
 {hyperref}
\usepackage{cite}
\hypersetup{%
pdfsubject={MSC 2000 Primary: 13P10; 90C57; Secondary: 65H10;12Y05; 90C27; 68R05},
pdfkeywords={order ideal polytope, border bases, Gröbner bases, combinatorial
optimization},
pdfauthor={Gábor Braun and Sebastian Pokutta}}

\numberwithin{equation}{section} 
\numberwithin{figure}{section} 
\theoremstyle{plain}
\newtheorem{thm}{Theorem}[section]
\newtheorem{corollary}[thm]{Corollary}
\newtheorem{lem}[thm]{Lemma}
\newtheorem{prop}[thm]{Proposition}

\theoremstyle{definition}
\newtheorem{defn}[thm]{Definition}
\newtheorem{example}[thm]{Example}

\theoremstyle{remark}
\newtheorem{rem}[thm]{Remark}

\theoremstyle{plain}
\newtheorem{inneralgorithm}[thm]{Algorithm}

\newenvironment{algorithm}[4]{%
  \begin{inneralgorithm}[{#1\textemdash#2}]\mbox{}%
      \begin{description}[itemsep=0pt]%
      \item[Input] #3
      \item[Output] #4
      \end{description}
    \begin{enumerate}[itemsep=0pt, topsep=0pt]
    }%
    {%
    \end{enumerate}%
  \end{inneralgorithm}}

\newtheoremstyle{problem}
  {}
  {}
  {\itshape}
  {}
  {\scshape}
  {:}
  { }
  {\thmnote{#3}}
\theoremstyle{problem}
\newtheorem{problem}{Problem}

\newcommand{\problemref}[1]{\textsc{\nameref{#1}}}

\newcommand*{\rk}[1]{\operatorname{rk}(#1)}

\DeclareMathOperator*{\conv}{conv}

\newcommand{\Z}{\mathbb{Z}}
\newcommand{\Q}{\mathbb{Q}}
\newcommand{\R}{\mathbb{R}}
\newcommand{\Pclass}{\textrm{P}\xspace}
\newcommand{\NP}{\textrm{NP}\xspace}
\newcommand{\coNP}{\textrm{coNP}\xspace}
\newcommand{\N}{\mathbb{N}}

\newcommand{\bT}{\mathbb{T}}
\DeclareMathOperator{\supp}{supp}
\newcommand{\oi}{\mathcal{O}}
\DeclareMathOperator{\LT}{LT}
\DeclareMathOperator{\LC}{LC}
\DeclareMathOperator{\LF}{LF}
\newcommand*{\textprogram}[1]{\textnormal{\texttt{#1}}}
\newcommand*{\textalgorithm}[1]{\textnormal{\texttt{#1}}}
\newcommand{\bbasis}{\ensuremath{\operatorname{%
      \hyperref[alg:borderBasisNew]{\textalgorithm{BBasis}}}}}
\newcommand{\lstabspan}{\ensuremath{\operatorname{%
      \hyperref[alg:stableSpan]{\textalgorithm{LStabSpan}}}}}
\newcommand{\basisT}{\ensuremath{\operatorname{%
      \hyperref[alg:basisT]{\textalgorithm{BasisTransformation}}}}}
\newcommand{\gausEl}{\ensuremath{\operatorname{%
      \hyperref[alg:gaussEl]{\textalgorithm{GaussEl}}}}}
\newcommand{\pring}{K[\mathbb{X}]}
\newcommand{\vgen}[1]{{\left\langle #1\right\rangle}\sb{K}}
\newcommand{\igen}[1]{{\left\langle #1\right\rangle}\sb{\pring}}

\DeclareMathOperator*{\argmax}{argmax}
\DeclarePairedDelimiter{\size}{|}{|}

\DeclareMathOperator{\xc}{xc}
\DeclareMathOperator{\OIP}{OIP}
\DeclareMathOperator{\CORR}{COR}

\newcommand{\allOne}{\mathbb{1}}

\begin{document}
\renewcommand{\sectionautorefname}{Section}
\renewcommand{\subsectionautorefname}{Section}

\title{A polyhedral characterization of border bases}

\author{Gábor Braun}
\address{ISyE,
  Georgia Institute of Technology,
  Atlanta, GA 30332, USA.}
\email{gabor.braun@isye.gatech.edu}

\author{Sebastian Pokutta}
\thanks{Research partially supported by German Research Foundation
  (DFG) funded SFB 805, Hungarian Scientific Research Fund,
  grant No. K 67928, NSF grants CMMI-1333789 and CCF-1415496}

\address{ISyE,
  Georgia Institute of Technology,
  Atlanta, GA 30332, USA.}
\email{sebastian.pokutta@isye.gatech.edu}

\subjclass[2000]{Primary: 13P10; 90C57; secondary: 65H10;12Y05; 90C27; 68R05}

\date{\today}

\keywords{order ideal polytope, border bases, Gröbner bases, combinatorial
optimization}

\begin{abstract}
  Border bases arise as a canonical generalization of Gröbner
  bases, using order ideals instead of term orderings.
  We provide a polyhedral characterization of all order ideals
  (and hence all border bases)
  that are supported by a zero-dimensional
  ideal: order ideals that support a border basis correspond
  one-to-one to integral points of the order ideal polytope.  In
  particular, we establish a crucial connection between the ideal and
  its combinatorial structure.  Based on this characterization we also
  provide an adaptation of the border basis algorithm of Kehrein and
  Kreuzer \cite{kehrein2006cbb} to allow for computing border bases
  for arbitrary order ideals,
  given implicitly via maximizing a preference on monomials
  (variable selection problem),
  \emph{independent} of term
  orderings.
  The algorithm requires the same size of resources
  as the border basis algorithm except for some minor overhead. We
  also show that the underlying variable selection problem of
  finding an order ideal that supports a border basis is NP-hard and
  that any linear description of the associated convex hull of all order ideals requires a
  superpolynomial number of inequalities.
  \end{abstract}
\maketitle

\section{Introduction}

In many different disciplines and real-world applications one is faced
with solving systems of polynomial equations. Often this is simply due
to a physical or dynamical system having a \emph{natural}
representation as a system of polynomial equations,
but equally often it
is due to the sheer expressive power of polynomial systems of
equations that allow for easy reformulation. To give an example of the
latter, an inequality \(ax \leq b\)
with \(a \in \R^n\)
and \(b\in \R\)
can be expressed via a single polynomial equation: \(ax + u^{2} = b\).
A slightly more involved example is that of the feasible region of a
binary program \(\{x \mid Ax \leq b, x \in \{0,1\}^n\}\),
which can be captured via rewriting each individual inequality
as before,
and adding quadratic polynomials \(x_i^2 - x_i = 0\)
for each coordinate \(i = 1, \dots, n\) of \(x\).
As a consequence,
there is a huge need to \emph{computationally} model,
understand, manipulate, and extract the solution set of systems of
polynomial equations.

A key insight in (computational) commutative algebra is that one can
choose a smart ordering on the monomials and compute a special set of
generators of the ideal generated by the system of equations
that makes many operations easy, and provides a structural insight into
the system.
One such special set of generators is
a Gröbner basis. By now, Gröbner bases are fundamental and standard
tools in commutative algebra to actually perform important operations
on ideals such as intersection, membership test, elimination,
projection, and many more.
Border bases arise as a
natural generalization of Gröbner bases that can be computed for
\emph{zero-dimensional ideals}, i.e., the associated factor space is a
\emph{finite-dimensional} vector space (see \autoref{sec:Preliminaries} for
details). While this might seem to be a severe restriction, for many
applications it is sufficient.
Roughly speaking, whenever the solution set is
finite, we are dealing with a zero-dimensional ideal. For example,
systems of polynomial equations with solutions restricted to a finite
set of points are captured by
zero-dimesional ideals. 

The advantage of border bases over Gröbner bases partly arises from
the iterative generation of linear syzygies, inherent in the border
basis algorithm, which allows for successively approximating the basis
of interest degree-by-degree, which leads to an implied notion of
approximability.  Moreover, many border basis algorithms (and also
Gröbner basis algorithms) are essentially linear algebra algorithms,
allowing for fast computation. However arguably the most important
difference between a Gröbner basis and a border basis is that the
former is computed with respect to a \emph{term ordering} (the
aforementioned ordering on the monomials) and the latter is computed
with respect to a so-called \emph{order ideal}, which for now can be
thought of as a generalization of a term ordering
(and not an ideal in
the usual sense).  An order ideal that \emph{supports} a border basis,
i.e., for which one can compute a border basis is called
\emph{admissible}. It is well known that every reduced Gröbner basis
can be extended to a border basis (see \cite[p.\
281ff]{kehrein2006cbb}), i.e., every term ordering
gives rise to an admissible order ideal. At the same time, not every
border basis is an extension of a Gröbner basis, since the former form a
strictly larger set giving potential extra freedom in modeling
solution sets to polynomial systems. For a given zero-dimensional
ideal \(I\)
of a polynomial ring \(R\),
the size of an order ideal that supports a border basis is
predetermined as the dimension of the vector space \(R / I\).
This is only a necessary condition though and not every order ideal
$\oi$ of size equal to the dimension of \(R / I\)
supports a border basis.  An example illustrating these two cases is
presented in \cite[Example 6]{kehrein2006cbb}. Finally, border bases
deform more smoothly in the input \cite{kreuzer2008dbb} (see also
\cite{robbiano2009border}), which is particularly helpful when the
coefficients arise from measurement data \cite{HKPP2009,
  abbott2008sbb}, e.g., that is why border bases are used in the
context of total-least-squares polynomial regression (see
\cite{HKPP2009} for details).

\subsection{Arbitrary order ideals}
\label{sec:arbitr-order-ideals}

While the above examples highlight the advantages of border bases for
certain types of computations, so far we have not answered a key
question: why it is desirable and important to be able to
compute border bases with respect to general order ideals (i.e., those
that do not necessarily stem from a term ordering).
We offer three
different perspectives.

First, it is
desirable to obtain a complete characterization of all border bases
supporting order ideals. In particular, our characterization can be
used to rule out certain types of order ideals \emph{and}
provides a \emph{proof} (a dual certificate) for their non-existence
via the associated violated inequality: in complexity-theoretic terms,
we provide certificates for the non-membership problem. Second,
choosing a different order ideal might significantly reduce
computational time to obtain a border basis. This aspect is well known
and often exploited in the context of
Gröbner bases. 

However, arguably the most important aspect from the perspective of
actually solving polynomial system is \emph{variable selection}. Often
the polynomial system of interest stems from e.g., a physical system
and the variables correspond to actual physical quantities and hence
have explanatory power. Now, it can be very desirable to obtain a
polynomial description of the solution set of the system using
\emph{specific monomials or variables} to allow for actual real-world
interpretation of the solutions. At the same time the ideal structure
might preclude a description with all desired variables or monomials
contained in the order ideal. We end up with an optimization problem
of finding an order ideal, which matches our preferences as well as
possible. Optimization problems of this type are referred to as
\emph{variable selection or feature selection} (see e.g,
\cite{james2013introduction}) and are ubiquitous in many data related
disciplines, such as e.g., statistics, machine learning, and more
broadly data analytics, where we effectively seek an explanation of a
phenomenon in specified explanatory variables. In our context this
naturally leads to the \problemref{pr:MaxOrderIdeal} problem,
where we specify weights for each monomial and we search a
maximal weight order ideal supporting a border basis of a given
zero-dimensional ideal \(I\).

Order ideals and determining those with maximum weight do appear in a
very natural way in combinatorial optimization as the so called
\emph{maximum weight closure problem} (a simplified version
of our \problemref{pr:MaxOrderIdeal} problem, cf.\ e.g.,
\cite{picard1976maximal}) and they have a variety of applications,
e.g., in open-pit mining where any feasible production plan is indeed
an order ideal; we refer the interested reader to
\cite{hochbaum2000performance} for an overview.
Another example is
the approximate vanishing ideal algorithm in \cite{HKPP2009},
which computes
a polynomial description of an approximate vanishing ideal of a given
set of (noisy) points. Effectively, a total-least-squares optimization
problem is solved here and explanatory variables can come from an
order ideal. If the points stem from actual (physical) measurements,
choosing the variables in the order ideal can help recover
important physical relations.
Other, more involved
applications might arise, e.g., in computational biology where the
structure of a boolean network is inferred from the Gröbner fan.

It has been an open question to characterize the
admissible order ideals of a zero-dimensional ideal. We provide
a polyhedral relaxation of all admissible order ideals of a given
zero-dimensional ideal that support a border basis. This is the best
we can hope for given that the separation
problem for the polytope is NP-hard
as we will see in
\autoref{sec:complexity-finding-ord-ideal}. Moreover, we will also
establish that in general the convex hull of all order ideals of a
given ideal can require a superpolynomial number of inequalities in
\emph{any} linear programming formulation, i.e., that the polytope
arising from the convex hull has superpolynomial extension
complexity. Many of the
results that we present later are the border bases analogs of
their counterpart
for Gröbner bases in \cite{onn1999cutting}.

\subsection{Computing border bases} The border basis algorithm in \cite{kehrein2006cbb}, which
is a specification of Mourrain's generic algorithm
\cite{mourrain1999ncn}, allows for computing border bases of
zero-dimensional ideals for order ideals supported by a
\emph{degree-compatible} term ordering.
However, this border basis algorithm does
not allow for the computation of a border basis
for more general order ideals (in fact it requires a degree-compatible
term ordering). The
alternative algorithm presented in \cite[Proposition
5]{kehrein2006cbb} which can \emph{potentially} compute arbitrary
border bases requires the \emph{a priori} knowledge of the order ideal
that might support a border basis. So while in principle the algorithm
can compute arbitrary border bases, the supporting order ideal has to
be part of the input.
Thus it does not characterize order ideals
for which a border basis does exist. Further, as pointed out in
\cite[p.\ 284]{kehrein2006cbb}, the basis transformation approach of
this algorithm is unsatisfactory as it significantly relies on Gröbner
basis computations. Another interesting approach for the computation
of normal forms that do not require degree-compatible term
orderings is \cite{mourrain2005generalized,mourrain2008stable},
however here a fixpoint scheme is required.

\subsection{Applications of border bases}
\label{sec:appl-bord-bases}

Surprisingly, it turns out that there are deep connections
to other mathematical disciplines and border bases represent the combinatorial
structure of the ideal under consideration in a canonical way. Although
the use of border basis as a concise framework is quite recent (see
e.g., \cite{kehrein2005cbb,kehrein2006cbb,kehrein14asv}),
the concept of border bases and in particular the border basis algorithm 
is rather old and has been reinvented
in different fields of mathematics
including computer algebra, discrete optimization, logic, and cryptography
under different names. In summary, border bases have been successfully
used for solving zero-dimensional systems of polynomial equations
(see, e.g., \cite{auzinger1988eac,moller1993sae,mourrain1999ncn}),
which in particular include those with 0/1 solutions and thus a
large variety of combinatorial problems.

\subsubsection*{Polynomial method}

Polynomial systems have been used in discrete mathematics and
combinatorial optimization to formulate combinatorial problems such as
the graph coloring problem, the stable set problem, and the matching
problem (we refer the interested reader to \cite{deLoera2009ecps}) as
well as to recognize graph properties \cite{de2010recognizing}.
This well-known method, which Alon referred to as the \emph{polynomial
  method} \cite{alon1999PolyMethod,alon1996PolyMethod} recently
regained strong interest and emphasizes the alternative view of
border bases algorithms in their various incarnations as proof systems which successively uncovers
hidden information by making it explicit.  In \cite[Section
2.3]{deLoera2009cpe} and \cite{deLoera2009ecps,deLoera2008hna,de2011computing}
infeasibility of certain combinatorial problems, e.g., 3-colorability
of graphs is established using Hilbert's Nullstellensatz and
the algorithm \textalgorithm{NulLA} is provided to establish
infeasibility by using a linear relaxation. The core of the algorithm
is identical to the \hyperref[alg:stableSpan]{$L$-stable span
  procedure} used in the border basis algorithm in \cite{kehrein2006cbb}, which intimately
links both procedures. The difference is of a technical but important nature:
whereas \textalgorithm{NulLA} establishes infeasibility, the 
border basis algorithm as presented in \cite{kehrein2006cbb} computes
the actual border bases of the ideal.  Another recent link between
border bases and the Sherali\textendash Adams closure \cite{SA} is
that the
Sherali\textendash Adams procedure can be understood as a weaker
version of the $L$-stable span procedure, see \cite{SP20092}.
Recently, border bases have
also been used to obtain a hierarchy of relaxations for polynomial
optimization problems \cite{bucero2014border}.

\subsubsection*{Border bases and cryptography}
Border bases have also been used to solve sparse quadratic systems of
equations thus giving rise to applications in cryptography in a
natural way. Such systems arise from crypto systems (such as e.g.,
AES, BES, HFE, DES, CTC variants) when rewriting the S-boxes as
polynomial equations. The celebrated XL, XSL, MutantXL attacks (see e.g.,
\cite{kreuzerAlgAttack,courtois2000eas})
are equivalent to the reformulation-linearization-technique (RLT) of Sherali and Adams
\cite{SA} and use a version of the Nullstellensatz to break
ciphers. 
Motivated by the success of the aforementioned methods,
border bases have also been used in cryptanalysis and coding theory,
see \cite{borgesquintana:ams}.

\subsubsection*{Border bases and numerical computations} Another core
application of border bases is the modeling of dynamic systems from
measured data (see e.g., \cite{HKPP2009, kreuzer2009sbb,
  abbott2008sbb}) where better numerical stability is
advantageous. The obtained solutions via border bases often provide a better
\emph{generalization}, i.e., explain new phenomena better,
than the respective Gröbner basis analog.

\subsection*{Our contribution}

Our contribution is the
following:

\subsubsection*{Polyhedral characterization of all border bases}
We provide a complete, polyhedral characterization of all border bases
of any zero-dimensional ideal $I$. We associate an \emph{order ideal polytope
$P$} to \(I\) whose integral points are in one-to-one correspondence with
order ideals supporting a border basis of $I$
(Theorem~\ref{thm:polytope=order-ideal}).
This explicitly establishes the link between the combinatorial structure
of the basis of the factor space and the structure of the ideal: whether
an order ideal supports a border basis is solely determined
by the combinatorial structure of the order ideal polytope. A related
result for Gröbner bases of the vanishing ideal of generic points was established in
\cite{onn1999cutting}, where it was shown that distinct reduced
Gröbner bases of the vanishing ideal are in bijection with the
vertices of the corner cut polyhedron. 

\subsubsection*{Computing maximum weight order ideals}
We will show that computing a maximum weight order ideal supporting a
border basis, i.e., solving the variable selection problem, is
\NP-hard in general (Theorem~\ref{thm:MaxOrderIdeal-NPhard}).
This is surprising as we merely
ask for a basis transformation. In particular, the \NP-hardness does
not stem from the hardness of computing the
\hyperref[alg:stableSpan]{\(L\)-stabilized span}, as the problem
remains \NP-hard, even in cases where the \(L\)-stabilized span is
small enough to be determined efficiently as shown in our reduction in
\autoref{sec:np-hard-constraints}.
In particular, unless \(\NP = \coNP\),
the convex hull of characteristic vectors of order ideals cannot have
an efficient linear programming formulation. In
Section~\ref{sec:extens-compl-order}, we complement this result and show
that there exists zero-dimensional ideals so that 
any linear programming formulation capturing their admissible order
ideals requires a subexponential
number of inequalities, irrespective of \(\NP\) vs. \(\coNP\); in the
language of extended formulations, we show that the convex hull of
admissible order ideals has subexponential extension complexity.
We discuss implications
of this in \autoref{sec:compl-order-ideal}.

\subsubsection*{Computing arbitrary border bases} We extend the
border basis algorithm in \cite{kehrein2006cbb}
to compute border bases for arbitrary order ideals using the order
ideal polytope (Algorithm~\ref{alg:borderBasisNew}),
where the order ideals are given implicitly by a preference vector.
(Note that every admissible order ideal can be obtained by choosing a
suitable preference vector.)
We would like to point out that algorithms
for general bases of quotient spaces have been proposed in
\cite{mourrain2005generalized,mourrain2008stable}). However these
algorithms are markedly different in relying on some fix point scheme,
so that the advantage of the degree-by-degree iterations are lost. We
refer the reader to the discussion in \cite{kaspar2013computing}.

\subsubsection*{Computational feasibility}

We provide computational tests that demonstrate the feasibility of our
method. Having the order ideal polytope available for a
zero-dimensional ideal $I$, it is possible to examine the structure of
the ideal based on its border bases. To demonstrate feasibility, we
consider the straightforward application of \emph{counting} the number
of border bases for a zero-dimensional ideal $I$, which we present as
an example in \autoref{sec:Computational-results} for counting
degree-compatible order ideals.

\subsection*{Subsequent work}
\label{sec:subsequent-works}

Following our work \cite{BP2009,BP20092}
several related results have been obtained. In \cite{ananth2011border}
it was shown that the border basis detection problem is NP-hard (see
also \cite{ananth2012complexity}). Moreover, an algorithm for computing
border bases without term orderings has been given in
\cite{kaspar2013computing}.

\subsection*{Outline}

We start with the necessary preliminaries in
\autoref{sec:Preliminaries} and recall the border basis algorithm from
\cite{kehrein2006cbb} in \autoref{sub:classicBB}.  In
\autoref{sec:The-order-ideal} we introduce the order ideal polytope
and establish the one-to-one correspondence between the integral
points of this polytope and border bases. We also derive an equivalent
characterization that is better suited for actual computations.  In
\autoref{sec:Computing-border-bases-deg-comp} we then use the results
from \autoref{sec:The-order-ideal} to obtain the generalized border
basis algorithm for arbitrary order ideals.  We then study the
complexity of the \problemref{pr:MaxOrderIdeal} problem establishing
\NP-hardness and a subexponential lower bound
on its polyhedral complexity in
\autoref{sec:complexity-finding-ord-ideal}.  We conclude with
computational results in \autoref{sec:Computational-results} and
with a summary in \autoref{sec:Concluding-remarks}.

\section{Preliminaries}
\label{sec:Preliminaries}

We consider a polynomial ring $\pring$ over the field $K$ with
variables
$\mathbb{X}=\{x_{1},\dots,x_{n}\}$.
Let $\bT^{n}\coloneqq \{\prod_{i} x_{i}^{k_{i}} \mid k_{i} \in\ N\}$
be the \emph{set of terms}, i.e., the set of all monomials.
Recall that the (total) degree of
a monomial \(m = \prod_{i} x_{i}^{k_{i}}\)
is \(\deg m = \sum_{i} k_{i}\).
For any $d\in\N$ we let $\bT_{\leq d}^{n}\coloneqq
\{m \in \bT^{n} \mid \deg m \leq d\}$
be the set of monomials of total degree at most $d$. Sometimes we will
refer to a subset of monomials \(L\) as the \emph{computational
  universe}, to which the actual computation is confined. For a polynomial
$p = \sum_{m \in \bT^{n}} a_{m} m \in \pring$ we define the
\emph{support of $p$} to be $\supp(p)\coloneqq \{m \in \bT^{n} \mid
a_{m} \neq 0\}$
and similarly, for a set of polynomials $P\subseteq\pring$ we define
the \emph{support of $P$} to be $\supp(P)\coloneqq \bigcup_{p\in P}\supp(p)$.
Given a (total) ordering \(\sigma\) on \(\bT^{n}\),
the \emph{leading term} \(\LT_{\sigma}(p) \coloneqq m\)
of the polynomial $p$
is the largest element \(m\) of \(\supp(p)\)
in the ordering \(\sigma\),
and the \emph{leading coefficient} $\LC_{\sigma}(p) = a_{m}$ of $p$
is the coefficient of $\LT_{\sigma}(p)$.
We drop the index \(\sigma\) if
the ordering is clear from the context.
Recall that a \emph{term ordering} is a total ordering \(\sigma\)
on \(\bT^{n}\) with \(m_{1} \leq m_{2}\) for all pair of monomials
\(m_{1}\), \(m_{2}\) with \(m_{1} \mid m_{2}\).
Monomial orderings are used for Gröbner basis computations,
but here we allow more general orderings.

The \emph{leading form} $\LF(p)$ of a polynomial
$p = \sum_{m \in \bT^{n}} a_{m} m \in \pring$
is defined to be $\LF(p) = \sum_{m : \deg m = \deg p} a_{m} m$,
i.e., we single out the part with maximum degree.
(The leading form does not depend on an ordering.)
Both $\LF$ and $\LT$ generalize to sets in the obvious way, i.e.,
for a set of polynomials $P$ we define $\LF(P)\coloneqq \{\LF(p)\mid p\in P\}$
and $\LT(P)\coloneqq \{\LT(p)\mid p\in P\}$.

In the following we
will frequently switch between considering polynomials $M$,
the generated ideal, and the generated vector space whose coordinates
are indexed by the monomials in the support of $M$. We denote the ideal
generated by $M$ as $\igen{M}$ and the vector space generated by
$M$ as $\vgen{M}$. For $n\in\N$ we define $[n] \coloneqq \{1,\dots,n\}$.
All other notation is standard as to be found in \cite{cox2007iva,kreuzer2000cca};
we have chosen the border basis specific notation to be similar to the one
in \cite{kehrein2006cbb}; see also \cite{kreuzer2005cca} for a
broader exposition. 

\subsection{Order ideals}
\label{sec:order-ideals}

Central to our discussion will be the notion of an order ideal,
which is \emph{not} an ideal, but a set of monomials closed under
taking (monomial) factors:

\begin{defn}
Let $\oi$ be a finite subset of $\bT^{n}$. If for all $m \in \oi$
and $m' \in \bT^{n}$ such that $m' \mid m$ we have $m' \in \oi$, i.e.,
$\oi$ is closed under factors, then we call $\oi$ an \emph{order
ideal}. Furthermore, the \emph{border $\partial\oi$} of a non-empty
order ideal $\oi$ is the set of monomials
$\partial\oi \coloneqq \{x_{j} m \mid j \in [n],
m \in \oi\} \setminus \oi$.
As an exception,
we set $\partial\emptyset\coloneqq \{1\}$ for the empty order ideal.
\end{defn}
Examples of order ideals are \(\{1\}\),
\(\{1, x_{1}, \dots, x_{1}^{k}\}\), and
\(\{1, x_{1}, x_{1}^{2}, x_{1} x_{2}, x_{2}\}\). Recall that an ideal $I\subseteq\pring$ is \emph{zero-dimensional}, if and only if $\pring/I$
is finite dimensional. The $\oi$-border basis of a zero-dimensional
ideal $I$ is a special set of polynomials:
\begin{defn}\label{def:bbasis}
Let $\oi = \{m_{1}, \dots, m_{\mu}\}$ be an order ideal
with border $\partial\oi=\{b_{1},\dots,b_{\nu}\}$.
Further let $I\subseteq\pring$ be a zero-dimensional ideal, and
\(\mathcal{G} = \{g_{1},\dots,g_{\nu}\} \subseteq I\)
be a (finite) set of polynomials.
Then the set $\mathcal{G}$ is an \emph{$\oi$-border
basis} of \(I\) if:
\begin{enumerate}
\item \label{enu:bbform} the polynomials in $\mathcal{G}$ have the
  form $g_{j}=b_{j}-\sum_{i=1}^{\mu}\alpha_{ij} m_{i}$
for $j\in[\nu]$ and $\alpha_{ij}\in K$;
\item $\pring=I\oplus\vgen{\oi}$ as vector spaces.
\end{enumerate}
If there exists an $\oi$-border basis of $I$ then the order
ideal $\oi$ \emph{supports a border basis} of $I$,
equivalently, $\oi$ is \emph{admissible} for $I$.
Let \(\Lambda(I)\) denote
the set of admissible order ideals of \(I\).
\end{defn}
Note that any border basis \(\mathcal{G}\) of an ideal \(I\)
is actually generating \(I\) as an ideal, i.e.
\(\igen{\mathcal{G}}=I\).
A proof of this fact can be found in \cite[Proposition
4.3.2]{kehrein14asv}; for the sake of completeness we provide an alternative proof here.
Let us consider
the subspace \(\igen{\mathcal{G}} + \vgen{\oi}\)
spanned by the ideal \(\igen{\mathcal{G}}\)
generated by \(\mathcal{G}\)
and the order ideal \(\oi\).
This subspace is closed under
multiplication by the \(x_{i}\), and hence it is an ideal.
As it contains \(1\)
(being contained in either \(\mathcal{G}\) or \(\oi\)),
it must be the whole ring, and hence
\(I = \left( \igen{\mathcal{G}} + \vgen{\oi} \right) \cap I
= \igen{\mathcal{G}} + \left( \vgen{\oi} \cap I \right)
= \igen{\mathcal{G}}\)
by the modular law.
Recall that the modular law states \((A + B) \cap C = A + (B \cap C)\)
for all subspaces \(A\), \(B\), \(C\) of a vector space with
\(A \subseteq C\).

In particular, an order ideal
$\oi$ supports an $\oi$-border basis of $I$
if and only if
$\pring = I \oplus \vgen{\oi}$.
Moreover, for any given order ideal
$\oi$ and ideal $I$ the $\oi$-border basis of $I$ is unique as
$b_{j}$ has a unique representation in $\pring=I\oplus\vgen{\oi}$
for all $j\in[\nu]$. Furthermore, as $\pring=I\oplus\vgen{\oi}$
it follows that $\lvert \oi \rvert = \dim \vgen{\oi} $ is invariant for all
choices of $\oi$.
The requirement for $I$ being zero-dimensional is necessary
to ensure finiteness of the order ideal \(\oi\) and its border $\partial\oi$.

\begin{example}[Order ideals from Gröbner bases]
  \label{ex:Groebner}
  A common way to obtain an admissible order ideal \(\oi\)
  for a zero-dimensional ideal \(I\)
  is to compute the Gröbner basis
  \(\mathcal{G}=\{g_{1},\dots,g_{\nu}\}\subseteq\pring\)
  of \(I\)
  with an arbitrary term ordering \(\sigma\),
  and let
  \[\oi \coloneqq \{ m \in \bT \mid
  \forall j \colon \LT_{\sigma}(g_{j}) \nmid m\}\]
  consists of all monomials \(m\) not divisible by
  any leading term in the Gröbner basis.

  As a concrete example, Tables~\ref{tab:running-ideal1}
  and~\ref{tab:running-ideal2}
  list every admissible order ideal of the following ideals
  over \(K[x_{1}, x_{2}]\),
  indicating a term ordering providing them if any
  (the ground field \(K\) can be any field):
  \begin{align}
    \label{eq:running-ideal1}
    I_{1} &\coloneqq
    \igen{x_{1} x_{2} - x_{1} - x_{2} + 1, x_{2}^{2} + x_{1} - 1}
    \\
    \label{eq:running-ideal2}
    I_{2} &\coloneqq
    \igen{x_{1}^{2} - x_{1} x_{2}, x_{2}^{2} - x_{1} x_{2},
      \bT^{2}_{=3}}
  \end{align}
  We enumerate the order ideals of \(I_{i}\) as
  \(\oi_{i,1}\), \(\oi_{i,2}\), \dots\
  with the first index
  identifying the ideal \(I_{i}\).
  \begin{table}
    \centering
    \begin{tabular}{|c|c|c|c|}
      \hline
      Term ordering & Gröbner basis & Admissible order ideal &
      Border basis \\
      \hhline{|=|=|=|=|}
      lex &
      \begin{tabular}{c}
        \(x_{1} + x_{2}^{2} - 1\),
        \\
        \(x_{2}^{3} - x_{2}^{2}\)
      \end{tabular}&
      \(\oi_{1,1} = \{1, x_{2}, x_{2}^{2}\}\)
      &
      \begin{tabular}{c}
        \(x_{1} + x_{2}^{2} - 1\),
        \\
        \(x_{1} x_{2} + x_{2}^{2} - x_{2}\),
        \\
        \(x_{1} x_{2}^{2}\),
        \(x_{2}^{3} - x_{2}^{2}\)
      \end{tabular}
      \\
      \hline
      degrevlex &
      \begin{tabular}{c}
        \(x_{1}^{2} - x_{1}\),
        \\
        \(x_{1} x_{2} - x_{1} - x_{2} + 1\),
        \\
        \(x_{2}^{2} + x_{1} - 1\)
      \end{tabular}&
      \(\oi_{1,2} = \{1, x_{1}, x_{2}\}\)
      &
      \begin{tabular}{c}
        \(x_{1} x_{2} - x_{1} - x_{2} + 1\),
        \\
        \(x_{1}^{2} - x_{1}\),
        \\
        \(x_{2}^{2} + x_{1} - 1\)
      \end{tabular}
      \\
      \hline
    \end{tabular}
    \caption{\label{tab:running-ideal1}%
      All admissible order ideals of \(I_{1} =
      \igen{x_{1} x_{2} - x_{1} - x_{2} + 1, x_{2}^{2} + x_{1} - 1}\)
      over \(K[x_{1}, x_{2}]\).
      Here \(x_{1} > x_{2}\) in both term orderings.}
  \end{table}
  \begin{table}
    \centering
    \begin{tabular}{|c|c|l|c|}
      \hline
      Term ordering & Gröbner basis & Admissible order ideal &
      Border basis\\
      \hhline{|=|=|=|=|}
      lex or degrevlex (\(x_{1} > x_{2}\))&
      \begin{tabular}{c}
        \(x_{1}^{2} - x_{1} x_{2}\),
        \\
        \(x_{1} x_{2} - x_{2}^{2}, x_{2}^{3}\)
      \end{tabular}&
      \(\oi_{2,1} = \{1, x_{1}, x_{2}, x_{2}^{2}\}\) &
      \begin{tabular}{c}
        \(x_{1}^{2} - x_{2}^{2}\), \\
        \(x_{1} x_{2} - x_{2}^{2}\), \\
        \(x_{1} x_{2}^{2}\),
        \(x_{2}^{3}\)
      \end{tabular}
      \\
      \hline
      lex or degrevlex (\(x_{2} > x_{1}\))&
      \begin{tabular}{c}
        \(x_{2}^{2} - x_{1} x_{2}\),
        \\
        \(x_{1} x_{2} - x_{1}^{2}, x_{1}^{3}\)
      \end{tabular}&
      \(\oi_{2,2} = \{1, x_{1}, x_{2}, x_{1}^{2}\}\)
      &
      \begin{tabular}{c}
        \(x_{2}^{2} - x_{1}^{2}\), \\
        \(x_{1} x_{2} - x_{1}^{2}\), \\
        \(x_{1}^{2} x_{2}\),
        \(x_{1}^{3}\)
      \end{tabular}
      \\
      \hline
      none & — & \(\oi_{2.3} = \{1, x_{1}, x_{2}, x_{1} x_{2}\}\)
      &
      \begin{tabular}{c}
        \(x_{1}^{2} - x_{1} x_{2}\), \\
        \(x_{2}^{2} - x_{1} x_{2}\), \\
        \(x_{1}^{2} x_{2}\),
        \(x_{1} x_{2}^{2}\)
      \end{tabular}
      \\
      \hline
    \end{tabular}
    \caption{\label{tab:running-ideal2}%
      All admissible order ideals of \(I_{2} =
      \igen{x_{1}^{2} - x_{1} x_{2}, x_{2}^{2} - x_{1} x_{2},
        \bT^{2}_{=3}}\)
      over the ring \(K[x_{1}, x_{2}]\).}
  \end{table}
  We leave it to the reader to verify that there is no further
  admissible order ideal for \(I_{1}\) and \(I_{2}\).
  (For \(I_{1}\), Figure~\ref{fig:proj-space-ideal1} showing linear
  dependence relations between monomials in
  \(K[x_{1}, x_{2}] / I_{1}\) should be helpful.)
  In the case of \(I_{2}\)
  the last admissible order ideal \(\{1, x_{1}, x_{2}, x_{1} x_{2}\}\)
  does not come from any term ordering, as we will show now.
  First note that this order ideal is indeed admissible:
  a basis of \(K[x_{1}, x_{2}] / I_{2}\)
  is given by the image of the admissible order ideal
  \(\oi_{2,2} = \{1, x_{1}, x_{2}, x_{1}^{2}\}\).
  As \(x_{1}^{2} - x_{1} x_{2} \in I_{2}\),
  the image of \(x_{1}^{2}\) is the same
  as the image of \(x_{1} x_{2}\),
  thus the image of
  \(\oi_{2.3} = \{1, x_{1}, x_{2}, x_{1} x_{2}\}\)
  coincides with that of \(\oi_{2,2}\),
  and hence it is a basis of \(K[x_{1}, x_{2}] / I_{2}\).
  In particular, \(K[x_{1}, x_{2}] = I_{2} \oplus \vgen{\oi_{2,3}}\)
  holds,
  showing the admissibility of \(\oi_{2,3}\).

  If \(\oi_{2,3}\) came from a Gröbner basis \(\mathcal{G}\)
  for a term ordering with \(x_{1} > x_{2}\),
  then \(\mathcal{G}\) would contain a \(g \in \mathcal{G}\)
  whose leading term divides the leading term \(x_{1} x_{2}\)
  of \(x_{1} x_{2} - x_{2}^{2} \in I\),
  and therefore \(\LT_{\sigma}(g) \in \oi_{2,3}\),
  a contradiction.
  The argument is similar for term orderings with \(x_{1} < x_{2}\).
\end{example}

Clearly,
as a vector space, every ideal $I$ has a degree filtration
$I=\bigcup_{i\in\N}I^{\leq i}$
where $I^{\leq i}\coloneqq \{ p\in I\mid\deg(p)\leq i \}$. For a set of
monomials $\oi$ we define
$\oi^{=i} \coloneqq \{m \in \oi \mid \deg(m) = i\}$,
and similarly
$\oi^{\leq i} \coloneqq \{m \in \oi \mid \deg(m) \leq i\}$.
In the following we will
also consider the special class of order ideals preserving the
degree filtration, which are called \emph{degree-compatible}:

\begin{defn}
\label{def:indepSizes}Let $I\subseteq\pring$ be a zero-dimensional
ideal and let $\oi\subseteq\bT^{n}$ be an order ideal. Then
$\oi$ is \emph{degree-compatible} (to $I$) if
\begin{equation}
  \label{eq:degree-compatible}
\size{\oi^{=i}} =
\size{\bT_{=i}^{n}} -\dim \frac{I^{\leq i}}{I^{\leq i-1}}
\end{equation}
for all $i\in\N$. 
\end{defn}
Thus, the $\oi$-border basis of a zero-dimensional ideal $I$ with
respect to any degree-compatible order ideal $\oi$ has a pre-determined
size for each degree $i\in\N$. Intuitively, the degree-compatible
order ideals are those that correspond to degree-compatible
orderings on the monomials. The important difference is that the orderings
do not have to be term orderings. The definition above only requires
\emph{local} compatibility with multiplication if $\oi$ is a degree-compatible order
ideal and thus downwardly closed, i.e., if $p,q$
are polynomials and $\deg(p)<\deg(q)$ then $p\leq q$.
Not all order ideals are degree-compatible as we will 
see in the following two examples.

\begin{example}[Degree-compatible order ideals]
  For finding degree-compatible ideals,
  particularly suitable term orderings are
  the \emph{degree-compatible} ones,
  like deglex and degrevlex,
  where \(m_{1} > m_{2}\) for all monomials \(m_{1}, m_{2}\)
  with \(\deg(m_{1}) > \deg (m_{2})\).
  Using the Gröbner basis \(\mathcal{G}\) of an ideal \(I\)
  under a degree-compatible term ordering,
  the low-degree parts \(I^{\leq i}\) of \(I\) can be easily
  determined
  using the elements of \(\mathcal{G}\) of degree at most \(i\),
  leading to the formula
  \begin{equation*}
    \dim \frac{I^{\leq i}}{I^{\leq i-1}} =
    \size{\{ m \in \bT^{n}_{=i} : \exists g \in \mathcal{G} \colon
    \LT_{\sigma}(g) \mid m\}},
  \end{equation*}
  from which \eqref{eq:degree-compatible} easily follows
  for the order ideal \(\oi\) coming from \(\mathcal{G}\),
  i.e., admissible order ideals coming from
  a degree-compatible term ordering are degree-compatible.

  Now it is easy to check that in Example~\ref{ex:Groebner},
  all the admissible order ideals are degree-compatible
  except
  \(\oi_{1,1} = \{1, x_{2}, x_{2}^{2}\}\)
  for \(I_{1} =
  \igen{x_{1} x_{2} - x_{2} - x_{1} + 1, x_{2}^{2} + x_{1} - 1}\).
\end{example}

\begin{example}[Generic ideal]\label{exa:generic}
  Let \(k\) and \(n\) be positive integers and
  let \({\{a_{ij}\}}_{i \in [n], j \in [k]}\) be
  algebraically independent real numbers over \(\Q\).
  Let \(I\) be the ideal of polynomials in the variables \(x_1,\dotsc,x_n\)
  which are zero on the points \((a_{1j}, \dotsc, a_{nj})\)
  for \(j \in [k]\), i.e., it is the vanishing ideal of those points.
  Thus, the ideal is zero-dimensional, and \(\pring / I\) has dimension \(k\).

  Every \(k\) distinct monomials form a complementary basis of \(I\),
  since they are linearly independent on the \(k\) points
  \((a_{1j}, \dotsc, a_{nj})\).
  An equivalent formulation of linear independence is that
  the determinant of the matrix formed by the values of
  the monomials on these points is non-zero.
  The determinant is indeed non-zero,
  as it is a non-trivial polynomial of the algebraic independent \(a_{ij}\)
  with integer coefficients.

  In particular, \emph{every} order ideal of size $k$ is an order ideal of \(I\).
  The degree-compatible order ideals are the ones
  where the monomials have the least possible degree,
  i.e., consisting of
  all monomials of degree less than \(l\) and
  in addition
  \(k - \binom{n+l-1}{l-1}\) monomials of degree \(l\),
  where \(l\) is the smallest non-negative integer satisfying
  \(k \leq \binom{n+l}{l}\), i.e.,
  there are at least \(k\) monomials of degree at most \(l\).
\end{example}

\subsection{Computing stable spans}
\label{sub:classicBB}

Without proofs,
we recall the underlying stable span computation of the border basis algorithm
in \cite{kehrein2006cbb} as it will serve as a basis for
our algorithm.  The interested reader is referred to
\cite{kehrein14asv,kehrein2005cbb} for a general introduction to
border bases and to \cite{kehrein2006cbb} in particular for an introduction
to the border basis algorithm.

The border basis algorithm in \cite{kehrein2006cbb}
calculates border bases of zero-dimensional ideals with
respect to an order ideal $\oi$ which is induced by a degree-compatible
term ordering $\sigma$ by successively generating a vector space approximation of the ideal. These approximations are generated via the following vector space \emph{neighborhood extensions}:

\begin{defn}
\label{def:lstabSpan}
(cf.\ \cite[Definition 7.1 and
the paragraph preceding Proposition 13]{kehrein2006cbb})
Let $V \subseteq \pring$ be a vector space.
We define the \emph{neighborhood
extension of $V$} to be
\[V^{+}\coloneqq V+Vx_{1}+\dots+Vx_{n}.\]

For a finite set \(W\) of polynomials,
its \emph{neighborhood extension} is
\[ W^{+} = W \cup W x_{1} \cup \dotsb \cup W x_{n}. \]
\end{defn}
Note that for a given set of polynomials $W$ such that $\vgen{W}=V$
we have $\vgen{W^{+}}=\vgen{W}^{+}=V^{+}$
as multiplication with $x_{i}$ is a $K$-linear map.
It thus suffices to perform the neighborhood extension on a set of
generators $W$ of $V$. 

Let $F$ be a finite set of polynomials and let $L\subseteq\bT^{n}$
be an order ideal, representing our computational universe.
We would like to compute the ideal generated by \(F\)
restricted to our universe, i.e.,  \(\igen{F} \cap \vgen{L}\).
We are mainly concerned with finite sets $L\subseteq\bT^{n}$.

Note that $F\cap \vgen{L} = \{f\in F\mid\supp(f)\subseteq L\}$,
i.e., $F\cap \vgen{L}$ contains only those polynomials that lie in the vector
space generated by $L$.
Clearly, $\vgen{F} \cap \vgen{\bT^n_{\leq d}} = \vgen{F}^{\leq d}$.
Using neighborhood extension we define:
\begin{defn}(Cf.\ \cite[Definition 10]{kehrein2006cbb})
Let $L$ be an order ideal and let $F$ be a finite
set of polynomials such that $\supp(F)\subseteq L$.
The set \(F\) is \emph{$L$-stabilized} if
$\vgen{F^{+}} \cap \vgen{L} = \vgen{F}$.
The \emph{$L$-stable span} \(F_L\) of \(F\) is
the smallest vector space \(V\) containing \(F\)
satisfying \(V^{+} \cap \vgen{L} = V\).
\end{defn}
The basic example of an \(L\)-stabilized set is
a set of generators for the intersection
\(I \cap L\) of an ideal \(I\) with \(L\),
but not all \(L\)-stabilized sets have necessary this form.
For example, \(L\) itself is \(L\)-stabilized
for any order ideal \(L\).
For \(L = \{1, x_{1}, x_{2}, x_{1}^{2} x_{2}^{2}, x_{2}^{3}\}\),
the set \(\{x_{1} + x_{2}, x_{1}^{2}, x_{2}^{2}, x_{2}^{3}\}\)
is \(L\)-stabilized,
but  \(\{x_{1} + x_{2}, x_{1}^{2}\}\) is not.

The following simple observation will be helpful later. 

\begin{rem}
\label{rem:stableForVR}
The \(L\)-stable span of a finite set \(F\) depends only on
the vector space \(\vgen{F}\) spanned by \(F\),
as \(\vgen{F^{+}} = \vgen{F}^{+}\).
\end{rem}

A straightforward construction of the \(L\)-stable span of \(F\) is
to inductively define the following increasing sequence of vector spaces:
\[
F_{0} \coloneqq \vgen{F} \quad\text{and}\quad
F_{k+1} \coloneqq F_{k}^{+} \cap \vgen{L}\text{ for \(k>0\)}.
\]
The union $\bigcup_{k\geq0} F_{k}$ is the $L$-stable span \(F_{L}\) of \(F\).

In the following we will explain how the $L$-stable span can be computed
explicitly for $L=\bT_{\leq d}^{n}$. We will use a modified version
of Gaussian elimination as a tool, which allows us to extend a given
basis $V$ with a set $W$ as described in the following:
\begin{lem}
\label{lem:borderGauss}\cite[Lemma 12]{kehrein2006cbb} Let $V=\{v_{1},\dots,v_{r}\}\subseteq \pring \setminus\{0\}$
be a finite set of polynomials such that $\LT(v_{i})\neq\LT(v_{j})$
whenever $i,j\in[r]$ with $i\neq j$ and $\LC(v_{i})=1$ for all
$i\in[r]$.  Further let $G$ be a finite set
of polynomials. Then Algorithm \ref{alg:gaussEl} computes a finite
set of polynomials $W\subseteq \pring$ with
\begin{enumerate}
\item 
  $\LC(w)=1$ for all $w\in W$,
\item
  $\LT(u_{1})\neq\LT(u_{2})$ for any
  distinct $u_{1},u_{2}\in V\cup W$, and
\item
  $\vgen{V\cup W}=\vgen{V\cup G}$.
\end{enumerate}
($V$, $W$ may be empty.)
\end{lem}
\begin{algorithm}{Gaussian Elimination for polynomials}{\gausEl}%
{$V$, $G$ finite set of polynomials, \(0 \notin V\)
  (as in Lemma \ref{lem:borderGauss}).}%
{$W\subseteq \pring$ finite set of polynomials
  (as in Lemma \ref{lem:borderGauss}).}
\label{alg:gaussEl}

\item Let $H\coloneqq G$ and $\eta\coloneqq 0$.
\item \label{enu:step2}If $H=\emptyset$ then return $W\coloneqq \{v_{r+1},\dots,v_{r+\eta}\}$
and stop.
\item Choose $f\in H$ and remove it from $H$. Let $i\coloneqq 1$.
\item
  If $f=0$ then go to step \ref{enu:step2}.
\item \label{enu:step4}
  If $i>r+\eta$ then
  put $\eta\coloneqq \eta+1$ and let $v_{r+\eta}\coloneqq f/\LC(f)$.
  Go to step \ref{enu:step2}.
\item If $\LT(f)=\LT(v_{i})$ then replace $f$ with $f - \LC(f) \cdot v_{i}$.
Set $i\coloneqq 1$ and go to step \ref{enu:step4}.
\item Set $i\coloneqq i+1$. Go to step \ref{enu:step4}.
\end{algorithm}

We can now compute the $L$-stable span using the Gaussian elimination
algorithm \ref{alg:gaussEl}:
\begin{lem}
\cite[Proposition 13]{kehrein2006cbb}\label{lem:stableSpan} Let
$L=\bT_{\leq d}^n$ and
$F \subseteq \pring$ be a finite set of polynomials supported on \(L\).
Then Algorithm \ref{alg:stableSpan} computes a vector space basis
$V$ of $F_{L}$ with pairwise different leading terms.\end{lem}
\begin{algorithm}{$L$-stable span computation}{\lstabspan}%
{$F$, $L$ as in Lemma \ref{lem:stableSpan}.}%
{$V$ as in Lemma \ref{lem:stableSpan}.}
\label{alg:stableSpan}
\item $V\coloneqq \gausEl(\emptyset,F)$.
\item \label{enu:step2StabSpan}$W'\coloneqq \gausEl(V,V^{+} \setminus V)$.
\item\label{enu:step3StabSpan}
  $W\coloneqq \{w\in W'\mid\supp(w) \subseteq L\} = \{ w \in W' \mid
  \deg(w) \leq d \}$.
\item If $W \neq \emptyset$ set $V \coloneqq V \cup W$ and go to step
  \ref{enu:step2StabSpan}.
\item Return \(V\).
\end{algorithm}

The rationale for computing a stable span approximation is due to
the following proposition
that serves as a criterion for testing whether an order ideal $\oi$ supports a border basis.

\begin{prop}
\label{pro:bbApproximation}
\cite[Proposition 16]{kehrein2006cbb}
Let $L$ be an order ideal.
Further let $\tilde{I}$ be an \(L\)-stabilized generating vector subspace of a zero-dimensional
ideal $I\subseteq\pring$, i.e., ${\tilde{I}}^{+}\cap \vgen{L}=\tilde{I}$
and $\igen{\tilde{I}}=I$. If $\oi$ is an
order ideal such that $\vgen{L}=\tilde{I}\oplus\vgen{\oi}$ and $\partial\oi\subseteq L$
then $\oi$ supports a border basis of $I$.
\end{prop}

We obtain the following corollary, which will be helpful later.

\begin{corollary}
  \label{cor:idealApprox}
  Let \(\tilde{I}\) be an \(\bT^{n}_{\leq d}\)-stabilized vector space
  satisfying \(\tilde{I} + \vgen{\bT^{n}_{\leq d-1}} = \vgen{\bT^{n}_{\leq d}}\).
  Then \(\igen{\tilde{I}} \cap \vgen{\bT^{n}_{\leq d}} = \tilde{I}\)
  and \(\left. \pring \middle/ \igen{\tilde{I}} \right. \cong
  \left. \vgen{\bT^{n}_{\leq d}} \middle/ \tilde{I} \right.\).
\begin{proof}
We apply Proposition~\ref{pro:bbApproximation} with the choice
\(L \coloneqq \bT^{n}_{\leq d}\),
\(I \coloneqq \igen{\tilde{I}}\)
and
\(\oi \coloneqq \bT^{n}_{\leq d} \setminus \LT(\tilde{I})\)
where the leading terms are
with respect to any degree-compatible ordering
(i.e., \(m_{1} < m_{2}\) whenever \(\deg m_{1} < \deg m_{2}\)).
Clearly,
\(\vgen{\bT^{n}_{\leq d}}= \tilde{I} \oplus \vgen{\oi}\).
The condition
\(\tilde{I} + \vgen{\bT^{n}_{\leq d-1}} = \vgen{\bT^{n}_{\leq d}}\)
ensures that
\(\oi\) consists of monomials of degree less than \(d\),
so \(\partial\oi \subseteq \bT^{n}_{\leq d}\).
Hence the proposition applies,
and we obtain \(\pring = I \oplus \vgen{\oi}\).
Together with
\(\vgen{\bT^{n}_{\leq d}}= \tilde{I} \oplus \vgen{\oi}\)
this gives
\(I \cap \vgen{\bT^{n}_{\leq d}} = \tilde{I}\),
and
\(\left. \pring \middle/ \igen{\tilde{I}} \right. \cong
  \vgen{\oi} \cong
  \left. \vgen{\bT^{n}_{\leq d}} \middle/ \tilde{I} \right.\).
\end{proof}
\end{corollary}

For a worst-case upper bound on $d$,
we will use the dimension of \(\pring / I\).
The necessary technical background is the following lemma.
\begin{lem}
  \label{lem:dimension}
  Let \(I\) be a zero-dimensional ideal of \(\pring\),
  and let \(d \coloneqq \dim \pring / I\).
  Then
  \begin{enumerate}
  \item \label{item:2}
    \(\left. I^{\leq d} \middle/ I^{\leq d-1} \right.
    \cong \vgen{\bT^n_{= d}}\) and
  \item \label{item:3}
    \(\left. \vgen{\bT^n_{\leq d-1}} \middle/
      I^{\leq d-1} \right.
    \cong \pring / I\).
  \end{enumerate}
\begin{proof}
Choose a degree-compatible term ordering.
The associated order ideal (as every order ideal of size \(d\))
contains monomials of degree less than \(d\).
This proves \(I + \vgen{\bT^n_{\leq d-1}} = \pring\),
from which the statements easily follow via the modular law
(\((A + B) \cap C = A + (B \cap C)\) for all subspaces
with \(A \subseteq C\)).
\end{proof}
\end{lem}

\section{The order ideal polytope}
\label{sec:The-order-ideal}

We will now introduce the \emph{order ideal polytope} $P(I)$
that characterizes all admissible order ideals, i.e., order ideals
supporting a border basis for a given zero-dimensional ideal $I$, in
an abstract fashion completely independent of
the stable span approximation. Its role will be crucial for the
\hyperref[alg:borderBasisNew]{later computation} of border bases for
  general order ideals. We will first focus on its properties and
  structure, then in \autoref{sec:Computing-border-bases}, we will
  consider the computational aspect.  We will show that the integral
  points of the order ideal polytope \(P(I)\)
  are in bijection with the admissible order
  ideals of a given zero-dimensional ideal
  $I$. In order to do so, we approach the problem from a polyhedral
  point of view to capture the intrinsic combinatorics for
  the admissibility condition
  $\pring = I \oplus \vgen{\oi}$ on the
  one hand and $\oi$ being an order ideal on the other hand. The role
  of the polyhedral description becomes prominent in
  \autoref{sec:Computing-border-bases} when the directness of the sum 
  $I \oplus \vgen{\oi}$ is rephrased in the language of matrices and
  vector space bases.

\subsection{Theoretical point of view}
\label{sec:theor-viewp}

We start with defining \emph{the order ideal polytope}
whose integral solutions are exactly the characteristic vectors
of order ideals admissible for a fixed zero dimensional ideal \(I\).
As we will see, the defining inequalities express various properties
of admissible order ideals.
\begin{defn}
\label{def:orderIdealPolytopeAlg}
Let \(I\) be a zero-dimensional ideal.
Its \emph{order ideal polytope} $P(I)$ is defined by the following
system of inequalities
with variables \(z_{m}\) for \(m \in \bT^n_{\leq d-1}\),
where $d \coloneqq \dim \left. \pring \middle/ I \right.$.
\begin{subequations}
  \label{eq:OrderIdealPolytope2}
  \begin{align}
  z_{m_{1}} &\geq z_{m_{2}}
  & \forall m_{1},m_{2}\in \bT^n_{\leq d-1} \colon m_{1}\mid m_{2}
  \label{eq:orderIdeal2}\\
  \sum_{m\in \bT^n_{\leq d-1}}z_{m} &= d
  \label{eq:degreeCompatibility2}\\
  \label{eq:steinitz2}
  \sum_{m\in U}z_{m} &\leq
    \dim \left. \vgen{U \cup I} \middle/ I \right. 
    & \forall U \subseteq \bT^n_{\leq d-1}\colon \lvert U \rvert =
    d
  \\
  z_{m} &\in [0,1] \quad &\forall m\in \bT^n_{\leq d-1} \label{eq:0-1}.
  \end{align}
\end{subequations}
\end{defn}
To obtain a finite dimensional polytope, we bounded the degree of the
monomials by $\dim \pring / I$ from above.
This bound is large enough to contain
all occurring monomials as we will see below. In a first step we relate the order ideal polytope
with admissible order ideals.  Recall that \(\Lambda(I)\) denotes the set of
admissible order ideals, i.e., order ideals
supporting a border basis of a zero-dimensional ideal
\(I\).

\begin{thm}
  \label{thm:polytope=order-ideal}
  Let \(I\) be a zero-dimensional ideal.
  There is an explicit bijection \(\xi\) between
  the set $\Lambda(I)$ of admissible order ideals of \(I\) and
  the set of integral points of the order ideal polytope
  \(P(I)\) of \(I\).
  The bijection is given by
  \begin{gather*}
    \xi \colon
    P(I) \cap \{0,1\}^{\bT^{n}_{\leq d-1}} \to \Lambda(I)
    \\
    \xi(z) = \oi(z) \coloneqq
    \{ m \in \bT^{n}_{\leq d-1} \mid z_{m}= 1\}.
  \end{gather*}
\begin{proof}
We show that the domain \(P(I) \cap \{0,1\}^{\bT^{n}_{\leq d-1}}\)
of \(\xi\) is exactly the set of characteristic vectors
of all order ideals
\(\oi \subseteq \bT^{n}_{\leq d-1}\)
admissible for \(I\).
It will immediately follow that \(\xi\) is a well-defined bijection
onto \(\Lambda(I)^{\leq d-1}\),
the set of admissible order ideals of \(I\) with all monomials
having degree less than \(d\).
Actually, this is the set of all admissible order ideals of \(I\),
as every admissible order ideal \(\oi\) of \(I\) has size
the dimension \(d\) of the factor \(\pring / I\),
and hence can only contain monomials up to degree \(d-1\).

Let \(z \in \{0,1\}^{\bT^{n}_{\leq d-1}}\) be a 0/1 vector.
It is the characteristic vector of the set
\(\oi(z) \coloneqq \{ m \in \bT^{n}_{\leq d-1} \mid z_{m}= 1\}\).
Recall that the set \(\oi(z)\) is an admissible order ideal of \(I\)
if and only if the following hold:
\begin{enumerate}
\item\label{item:order-ideal}
  \(\oi(z)\) is an order ideal, i.e.,
  \(m_{2} \in \oi(z)\) implies \(m_{1} \in \oi (z)\)
  for all monomials \(m_{1}\)  and \(m_{2}\)
  with \(m_{1} \mid m_{2}\).
\item\label{item:size-dimension}
  \(\size{\oi(z)} = d\).
\item\label{item:disjoint}
  The image of \(\oi(z)\) in \(\pring / I\)
  is linearly independent.
\end{enumerate}
The last two conditions together are clearly
an equivalent formulation of \(\pring = \vgen{\oi(z)} \oplus I\),
using that \(\oi(z)\) is a set of linearly independent elements
in \(\pring\).

Now we rewrite these conditions for the characteristic vector \(z\).
Condition~\ref{item:order-ideal} for fixed monomials
\(m_{1} \mid m_{2}\) is obviously equivalent to
\(z_{m_{1}} \geq z_{m_{2}}\).
Therefore Condition~\ref{item:order-ideal} is equivalent to
\eqref{eq:orderIdeal2}.
Similarly,
as \(\sum_{m \in \bT^{n}_{\leq d-1}} z_{m} = \size{\oi(z)}\),
Condition~\ref{item:size-dimension} is equivalent to
\eqref{eq:degreeCompatibility2}.

As of Condition~\ref{item:disjoint}, we first give a more complex but
equivalent formulation:
\begin{equation}
  \label{eq:disjoint}
    \left\lvert U \cap \oi(z) \right\rvert \leq
    \dim \left. \vgen{U \cup I} \middle/ I \right.,
    \qquad \text{for all \(U \subseteq \bT^{n}_{\leq d-1}\) with
      \(\size{U} = d\).}
\end{equation}
i.e., the size of \( U \cap \oi(z)\) is at most
the dimension of the vector space generated by the image of \(U\) in
the factor
\(\pring / I\).
This is obviously \emph{necessary} for
the image of \(\oi(z)\) to be linearly independent in the factor,
as then the image of \(U \cap \oi(z)\) is independent,
and contained in \(\left. \vgen{U \cup I} \middle/ I \right.\).
(For necessity, the size of \(U\) does not matter.)
For \emph{sufficiency} choose $U \coloneqq \oi(z)$,
showing that the image of \(\oi(z)\) spans a subspace of
\(\pring / I\) of size at least that of \(\oi(z)\),
i.e., that the image of \(\oi(z)\) is linearly independent.
Thus Condition~\ref{item:disjoint} is equivalent to
\eqref{eq:disjoint},
which is \eqref{eq:steinitz2} using \(\left\lvert U \cap \oi(z)
\right\rvert = \sum_{m\in U}z_{m}\).

All in all, a 0/1 vector \(z \in \{0,1\}^{\bT^{n}_{\leq d-1}}\)
is the characteristic vector of an order ideal admissible for \(I\)
if and only if it satisfies \eqref{eq:orderIdeal2},
\eqref{eq:degreeCompatibility2} and \eqref{eq:steinitz2}.
In other words,
\(P(I) \cap \{0,1\}^{\bT^{n}_{\leq d-1}}\)
is the set of characteristic vectors of order ideals
\(\oi \subseteq \bT^{n}_{\leq d-1}\) admissible to \(I\),
as claimed.
(The remaining inequalities \eqref{eq:0-1} of \(P(I)\)
are satisfied by all 0/1 vectors.)
\end{proof}
\end{thm}

\begin{example}[Order ideal polytope]
As an easy example we determine the order ideal polytope \(P(I_{2})\)
of the ideal
\(I_{2} = \igen{x_{1}^{2} - x_{1} x_{2}, x_{2}^{2} - x_{1} x_{2},
  \bT^{2}_{=3}}\)
from Example~\ref{ex:Groebner}.
First we derive several valid inequalities for \(P(I_{2})\)
in order to obtain a simple description.

Recall that \(K[x_{1}, x_{2}] / I_{2}\) has dimension \(d=4\),
so the coordinates of the polytope are indexed by monomials up to
degree \(3\).
As the ideal \(I\) contains \(\bT^{2}_{=3}\),
by \eqref{eq:steinitz2} applied to \(U = \bT^{2}_{=3}\), which
consists of exactly \(4\) monomials, 
\begin{equation}
  \label{eq:polytope2-degree3}
  \sum_{m \in \bT^{2}_{=3}} z_{m} \leq 0.
\end{equation}
Together with \(z_{m} \geq 0\) for all \(m\),
this implies
\begin{equation*}
  z_{m} = 0, \qquad \text{whenever } \deg m = 3.
\end{equation*}
Hence from now on we can omit variables indexed by degree-three
monomials as they are \(0\).

Now we apply \eqref{eq:steinitz2} again,
but this time for \(U = \bT^{2}_{=2} \cup \{x_{1}^{3}\}\)
(the role of the monomial \(x_{1}^{3}\) is simply to pad \(U\) 
ensuring that \(U\) has \(4\) elements) and derive
\begin{equation*}
  z_{x_{1}^{2}} + z_{x_{1} x_{2}} + z_{x_{2}^{2}} \leq 1.
\end{equation*}
Together with
\begin{equation*}
  z_{1} + z_{x_{1}} + z_{x_{2}} +
  z_{x_{1}^{2}} + z_{x_{1} x_{2}} + z_{x_{2}^{2}} = 4
\end{equation*}
by~\eqref{eq:degreeCompatibility2}
and \(z_{m} \leq 1\)
for all \(m \in \bT^{2}_{\leq 1}\) by \eqref{eq:0-1},
we obtain
\begin{align*}
  z_{1} = z_{x_{1}} = z_{x_{2}}
  &
  = 1,
  \\
  z_{x_{1}^{2}} + z_{x_{1} x_{2}} + z_{x_{2}^{2}}
  &
  = 1.
\end{align*}

All in all, the polytope \(P(I_{2})\) satisfies
the following inequalities:
\begin{align*}
  z_{m} &= 0, & \text{whenever } \deg m &= 3,
  \\
  z_{1} = z_{x_{1}} = z_{x_{2}}
  &
  = 1,
  \\
  z_{x_{1}^{2}} + z_{x_{1} x_{2}} + z_{x_{2}^{2}}
  &
  = 1,
  \\
  z_{x_{1}^{2}}, z_{x_{1} x_{2}}, z_{x_{2}^{2}}
  &
  \geq 0.
\end{align*}
This system defines a triangle with the following vertices and gives
rise to the order ideals \(\oi_{2,2} = \{1, x_1,x_2,x_{1}^{2}\}\),
\(\oi_{2,3} = \{1, x_1,x_2,x_{1}x_2\}\), and \(\oi_{2,1} = \{1,
x_1,x_2,x_{2}^{2}\}\): 
\begin{itemize}
\item
  \((z_{x_{1}^{2}} = 1, z_{x_{1} x_{2}} = 0, z_{x_{2}^{2}} = 0)\) 
  \hfill [characteristic vector of \(\oi_{2,2}\)]
\item
  \((z_{x_{1}^{2}} = 0, z_{x_{1} x_{2}} = 1, z_{x_{2}^{2}} = 0)\)
  \hfill [characteristic vector of \(\oi_{2,3}\)]
\item
  \((z_{x_{1}^{2}} = 0, z_{x_{1} x_{2}} = 0, z_{x_{2}^{2}} = 1)\)
  \hfill [characteristic vector of \(\oi_{2,1}\)]
\end{itemize}
We conclude that the \(P(I_{2})\) is the triangle with vertices
the characteristic vectors of the admissible order ideals of
\(I_{2}\), listed in Table~\ref{tab:running-ideal2}.
\end{example}
\begin{example}[A non-integral order ideal polytope]
  In contrast to \(P(I_{2})\) of the previous example,
  the order ideal polytope \(P(I_{1})\) of \(I_{1}\) from
  Example~\ref{ex:Groebner} is not the convex hull of
  the characteristic vectors of the admissible order ideals
  of \(I_{1}\) but a proper relaxation, i.e., \(P(I_{1})\) is more than the line segment
  of \(\oi_{1,1}\) and \(\oi_{1,2}\).
  A point of \(P(I_{1})\) lying outside this segment is
  the one with all its coordinates being \(1/2\):
  \begin{equation}
    z_{1} = z_{x_{1}} = z_{x_{2}}
    =
    z_{x_{1}^{2}} = z_{x_{1} x_{2}} = z_{x_{2}^{2}} = 1/2.
  \end{equation}
  Recall that \(K[x_{1}, x_{2}] / I_{1}\) has dimension \(3\).

  This point satisfies the system \eqref{eq:OrderIdealPolytope2},
  from which only \eqref{eq:steinitz2} requires explanation.
  For computing the right-hand side of \eqref{eq:steinitz2},
  the key is to determine the linear dependence relations between the
  monomials of degree at most \(2\) in the factor space
  \(K[x_{1}, x_{2}] / I_{1}\).
  These can be easily read off from
  Figure~\ref{fig:proj-space-ideal1}
  depicting
  the projective space of the factor \(K[x_{1}, x_{2}] / I_{1}\),
  which has dimension \(2\).
  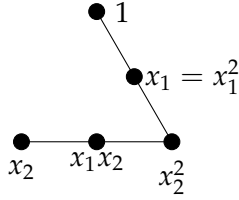
\begin{figure}
    \centering
    \begin{tikzpicture}[point/.style={shape=circle,fill,draw,
        minimum size=6pt, inner sep=0pt}]
      \draw (0,0) node[point,label=below:{\(x_{2}\)}]{}
      --
      node[point, label=below:{\(x_{1} x_{2}\)}] {}
      (2,0)
      node[point, label=below:{\(x_{2}^{2}\)}] {}
      --
      node[point, label=right:{\(x_{1} = x_{1}^{2}\)}] {}
      (60:2)
      node[point, label=right:{\(1\)}] {};
    \end{tikzpicture}
    \caption{\label{fig:proj-space-ideal1}%
      Linear dependence relations between
      low-degree monomials in the factor \(K[x_{1}, x_{2}] / I_{1}\)
      of dimension \(3\).
      The monomials are depicted in the projective space
      of the factor, to save a dimension.}
  \end{figure}
  To verify the figure,
  note that \(\oi_{1,1} = \{1, x_{2}, x_{2}^{2}\}\)
  is a basis of the factor,
  and hence forms a triangle in the projective space.
  As
  \(x_{1} + x_{2}^{2} - 1\),
  \(x_{1}^{2} - x_{1}\),
  \(x_{1} x_{2} + x_{2}^{2} - x_{2}\)
  are all elements of \(I_{1}\),
  we immediately see that \(x_{1} x_{2}\)
  is a third point on the line joining \(x_{2}\) and \(x_{2}^{2}\),
  and \(x_{1}\) is a third point on the line of
  \(1\) and \(x_{2}^{2}\),
  while \(x_{1}^{2}\) is the same point as \(x_{1}\).

  Now from Figure~\ref{fig:proj-space-ideal1}
  it is immediate that any subset \(U \subseteq \bT^{2}_{\leq 2}\)
  of size \(3\) has dimension at least \(2\)
  in \(K[x_{1}, x_{2}] / I_{1}\).
  Therefore the right-hand side of \eqref{eq:steinitz2},
  is at least \(2\),
  while the left-hand side is exactly \(3/2\),
  and therefore the inequality holds as claimed.
\begin{rem}
  Adding to \eqref{eq:OrderIdealPolytope2}
  the equality \(z_{1} = 1\) and requiring
  \eqref{eq:steinitz2} for all \(U \subseteq \bT^{n}_{\leq d-1}\)
  would still be insufficient to describe the convex hull
  of admissible order polytopes of \(I_{1}\):
  the following point still lies outside the convex hull
  while satisfying even the additional constraints:
  \begin{align}
    z_{0} &= 1, &
    z_{x_{1}} = z_{x_{2}}
    =
    z_{x_{1}^{2}} = z_{x_{1} x_{2}} = z_{x_{2}^{2}} &= 2/5.
  \end{align}
\end{rem}
\end{example}

\subsection{Computational point of view}
\label{sec:Computing-border-bases}

From a computational perspective,
the system \eqref{eq:OrderIdealPolytope2} defining the order ideal
polytope contains
dimensions \(\dim \left. \vgen{U \cup I} \middle/ I \right.\),
which are computationally challenging to determine.
Therefore in this subsection we provide a modified description
of the order ideal polytope, well suited for computations.

Let \(M \subseteq \pring\)
be a finite set of polynomials, and let \(M^{=i}\) denote
the set of polynomials in \(M\) with total degree \(i\).
We would like to have \(M\) to be a vector space basis
of \(I^{\leq d}\) reflecting the degree filtration of \(I\).
The following definition will be helpful. 

\begin{defn}
\label{def:L-canonical-form}
Let $M$ be a finite set of non-zero polynomials
of degree at most \(\ell\)
for some $\ell\in\N$.
Then \(M\) is in \emph{canonical form}
if the leading term of any element of \(M\) does not occur
in the other elements.
\end{defn}
Here the ordering on monomials can be any degree-compatible (total) ordering
(i.e., \(m_{1} < m_{2}\) for all monomials \(m_{1}\), \( m_{2}\)
with \(\deg m_{1} < \deg m_{2}\)),
and need not be a term ordering.
Clearly, any vector space basis can be brought into canonical form
via Gaussian elimination. The \emph{coefficient matrix}
 $A\in K^{M\times\bT_{\leq \ell}^{n}}$ of \(M\)
is the matrix where
the rows are indexed by the elements of \(M\),
and the columns are indexed by
the monomials of degree at most \(\ell\),
and the entries are the coefficients of the monomials in the elements of
\(M\).
In other words, \(A_{f,m} = a_{m}\)
for \(f = \sum_{m \in \bT^{n}_{\leq \ell}} a_{m} m \in M\). A visual
interpretation of a set \(M\) in canonical form
can be found in \autoref{fig:l-can-form}
using the coefficient matrix.

\begin{figure}
\[
A=\left(\begin{array}{ccc|c||ccc|c||c||ccc|c}
1 &  & 0 & \star &  &  &  & \star &  &  &  &  & \star\\
 & \ddots &  & \vdots &  & 0 &  & \vdots &  &  & 0 &  & \vdots\\
0 &  & 1 & \star &  &  &  & \star &  &  &  &  & \star\\
\hline  &  &  & 0 & 1 &  & 0 & \star &  &  &  &  & \star\\
 & 0 &  & \vdots &  & \ddots &  & \vdots &  &  & 0 &  & \vdots\\
 &  &  & 0 & 0 &  & 1 & \star &  &  &  &  & \star\\
\hline  &  &  &  &  &  &  &  & \ddots & & & & \\
\hline  &  &  & 0 &  &  &  & 0 &  & 1 &  & 0 & \star\\
 & 0 &  & \vdots &  & 0 &  & \vdots &  &  & \ddots &  & \vdots\\
 &  &  & 0 &  &  &  & 0 &  & 0 &  & 1 & \star\end{array}\right)\]

\caption{Coefficient matrix of a set of polynomials in canonical form.
  Double lines separate same-degree blocks of monomials.}
\label{fig:l-can-form}
\end{figure}

The following lemma summarizes the required properties of
a generating set $M$ of an ideal \(I\)
sufficient to describe the degree filtration
of \(I^{\leq d}\).

\begin{lem}
\label{lem:stableIdeal}
Let $M$ be in canonical form and \(\bT^{n}_{\leq d}\)-stabilized.
Further assume
\(\bT^{n}_{=d} \subseteq \vgen{M} + \vgen{\bT^{n}_{\leq d-1}}\).
Then the following hold for all $i\in[d]$:
\begin{enumerate}
\item\label{item:1}
  A basis for
  \(\left. \igen{M}^{\leq i} \middle/ \igen{M}^{\leq i-1} \right.\)
  is the image of \(M^{=i}\).
\item \(\igen{M}^{\leq i} = \vgen{\bigcup_{j \leq i} M^{=j}}\)
\item $\vgen{M^{=i}}^{<i} = 0$ and thus
  $\vgen{M^{=i}}^{<i} \subseteq \vgen{\bigcup_{0 \leq j \leq i-1} M^{=j}}$
\end{enumerate}
\begin{proof}
We first show that $\vgen{M^{=i}}^{<i}=0$ for all $i\in[d]$. Let
$i\in[d]$ be arbitrary and observe that each nonzero element $p\in M^{=i}$
has degree $i$. As $M$ is in canonical form, the polynomials
in $M^{=i}$ are interreduced (see the matrix in
Figure~\ref{fig:l-can-form} for Definition \ref{def:L-canonical-form})
and thus each nonzero element $p\in\vgen{M^{=i}}$ has
degree $i$.

By Corollary~\ref{cor:idealApprox},
we have \(\igen{M} \cap \vgen{L} = \vgen{M}\).
Hence \(\igen{M}^{\leq i} = \vgen{M}^{\leq i}\) for \(i \in [d]\).
Now the statements of the lemma are obvious consequences of \(M\)
being in canonical form.
\end{proof}
\end{lem}
The following lemma provides us a practical way to compute
the sizes of the degree components of degree-compatible order ideals,
which are the same for all order ideals of a given ideal.
\begin{lem}
\label{lem:stableRepresentation}
Let
$M$ be in canonical form and $\bT_{\leq d}^{n}$-stabilized.
Further let $\oi$ be an order ideal of $\igen{M}$,
and $d=\max_{m\in\partial\oi}\deg(m)$.
Let us assume \(\bT^{n}_{=d} \subseteq \vgen{M} + \bT^{n}_{\leq d-1}\).
Then $\oi$ is degree-compatible if and only if \[
\size{\oi^{=i}} = \size{\bT^{n}_{=i}} - \size{M^{=i}}\]
for every \(i \in [d]\).
\begin{proof}
In view of Definition \ref{def:indepSizes} it suffices to observe
that
$\left. I^{\leq i} \middle/ I^{\leq i-1} \right.$
has the image of \(M^{=i}\) as a basis
by Lemma~\ref{lem:stableIdeal} \ref{item:1}.
\end{proof}
\end{lem}

We are ready to provide a reformulation of
the definition of order ideal polytopes,
which is better suited for actual computations, partly as
it involves only direct matrix operations
via replacing dimensions with ranks of submatrices.
While \(d\) will still be the
dimension of \(\pring / \igen{M}\),
we do not require explicit \emph{a priori} knowledge,
but rather formulate alternative, sufficient conditions,
which are easier to verify by an algorithm.

\begin{lem}
\label{lem:polytope-comp}
Let $M$ be $\bT^{n}_{\leq d}$-stabilized and in canonical form
for some \(d \in \N\).
Suppose
\(\left. \vgen{M}^{\leq d} \middle/ \vgen{M}^{\leq d-1} \right. \cong
\vgen{\bT^n_{= d}}\) and \(d = \lvert \bT^{n}_{\leq d} \rvert - \lvert M \rvert\).
Then an alternative description
of the order ideal polytope
$P(\igen{M})$
of \(\igen{M}\)
is given by the system of inequalities
\begin{subequations}
  \label{eq:OrderIdealPolytope}
\begin{align}
  z_{m_{1}} &\geq z_{m_{2}} &
  \forall m_{1},m_{2}\in \bT_{\leq d-1}^{n}\colon m_{1}\mid m_{2}\label{eq:orderIdeal}\\
  \sum_{m\in \bT_{\leq d-1}^{n}}z_{m} &= d
  \label{eq:degreeCompatibility}\\
  \sum_{m\in U}z_{m} &\geq
  \lvert U \rvert -\rk{\tilde{U}}\label{eq:equiSteinitz}
  &
    \forall U\subseteq \bT_{\leq d-1}^{n} \colon
  \size{U}= \size{M^{\leq d-1}}
  \\
  0 \leq z_{m} &\leq 1 & \forall m\in \bT_{\leq d-1}^{n}.
\end{align}
In \eqref{eq:equiSteinitz}, the matrix
$\tilde{U}$ is the submatrix of
the coefficient matrix of $M^{\leq d-1}$
consisting of only the columns
indexed by monomials in $U$.
\end{subequations}
\begin{proof}
Let \(I \coloneqq \igen{M}\).
As \(\bT^{n}_{\leq d} = \vgen{M} + \bT_{\leq d-1}^{n}\) by assumption,
Lemma~\ref{lem:stableIdeal} provides
\(\igen{M} \cap \vgen{\bT^{n}_{\leq d}} = \vgen{M}\).
Moreover, \(\pring = \igen{M} + \bT_{\leq d-1}^{n}\),
hence via an application of the modular law
\(\left. (\igen{M} + \bT_{\leq d-1}^{n}) \middle/ \igen{M} \right.
= \left. \bT_{\leq d-1}^{n} \middle/
  (\igen{M} \cap \bT_{\leq d-1}^{n}) \right.\),
we obtain
\(\left. \pring \middle/ \igen{M} \right.
= \left. \bT_{\leq d-1}^{n} \middle/ \vgen{M^{\leq d-1}} \right.\).
In particular,
\(\dim \left. \pring \middle/ \igen{M} \right.
= \size{\bT_{\leq d-1}^{n}} - \size{M^{\leq d-1}} = d\).

As \(d = \dim \left. \pring \middle/ \igen{M} \right.\),
the only difference between the systems
\eqref{eq:OrderIdealPolytope2} and \eqref{eq:OrderIdealPolytope}
is that \eqref{eq:steinitz2} is replaced by
\eqref{eq:equiSteinitz}.
So we will show their equivalence modulo the other inequalities.

We start with \eqref{eq:equiSteinitz}
for a fixed $U \subseteq \bT_{\leq d-1}^{n}$, and make equivalent
transformations to it.
(For the following argument the size of \(U\) is irrelevant.)
Taking the difference with the equality
\eqref{eq:degreeCompatibility},
we obtain
\begin{equation}
  \label{eq:2}
  \sum_{m \in \bT^n_{\leq d-1} \setminus U} z_{m}
  \leq
  d - \size{U} + \rk{\tilde{U}}.
\end{equation}
Next we rewrite the right-hand side.
Recall that \(\tilde{U}\) is the submatrix obtained by restricting
to the columns corresponding to the monomials in \(U\),
i.e., the coefficient matrix of the image of \(M^{\leq d - 1}\)
in the factor
\(\left. \bT_{\leq d-1}^{n} \middle/
  \vgen{\bT_{\leq d-1}^{n} \setminus U} \right.\).
Therefore
\[
\rk{\tilde{U}} =
\dim \frac{\vgen{M^{\leq d-1} \cup (\bT_{\leq d-1}^{n} \setminus U)}}
{\vgen{\bT_{\leq d-1}^{n} \setminus U}}
=
\dim \vgen{M^{\leq d-1} \cup (\bT_{\leq d-1}^{n} \setminus U)}
- \size{\bT_{\leq d-1}^{n} \setminus U}.\]
Thus the right-hand side of \eqref{eq:2} becomes
\begin{equation}
  \label{eq:3}
 \begin{split}
  d - \size{U} + \rk{\tilde{U}}
  &
  =
  d - \size{U}
  +
  \dim \vgen{M^{\leq d-1} \cup (\bT_{\leq d-1}^{n} \setminus U)}
  - \size{\bT_{\leq d-1}^{n} \setminus U}
  \\
  &
  =
  d
  - \size{\bT_{\leq d-1}^{n}}
  + \dim \vgen{M^{\leq d-1} \cup (\bT_{\leq d-1}^{n} \setminus U)}
  \\
  &
  =
  - \dim \vgen{M^{\leq d-1}}
  + \dim \vgen{M^{\leq d-1} \cup (\bT_{\leq d-1}^{n} \setminus U)}
  \\
  &
  =
  \dim
  \frac{\vgen{M^{\leq d-1} \cup (\bT_{\leq d-1}^{n} \setminus U)}}
  {\vgen{M^{\leq d-1}}}
  =
  \dim \frac{\vgen{I \cup (\bT_{\leq d-1}^{n} \setminus U)}}{I}
  ,
 \end{split}
\end{equation}
where the last equality follows via the modular law
\[\frac{\vgen{M^{\leq d-1} \cup
    \left( \bT_{\leq d - 1}^{n} \setminus U \right)} + I}{I}
=
\frac{ \vgen{M^{\leq d-1} \cup
    \left( \bT_{\leq d - 1}^{n} \setminus U \right)}}
{\vgen{M^{\leq d-1} \cup (\bT_{\leq d - 1}^{n} \setminus U)} \cap I}
.\]
Therefore for a fixed \(U \subseteq \bT^{n}_{\leq d-1}\),
the inequality~\eqref{eq:equiSteinitz}
is equivalent to
\eqref{eq:steinitz2} with \(U\) replaced by
its complement \(\bT^{n}_{\leq d-1} \setminus U\).
The equivalence of \eqref{eq:equiSteinitz} and
\eqref{eq:steinitz2} stated for all subsets \(U\) follows,
noting that \(\size{M^{\leq d-1}} = \size{\bT^{n}_{\leq d-1}} - d\),
as shown at the beginning of the proof.
\end{proof}
\end{lem}

\section{\label{sec:Computing-border-bases-deg-comp}Computing border bases using the
order ideal polytope}

In the following we explain how Theorem \ref{thm:polytope=order-ideal}
can be used to actually compute border bases for general
order ideals. We cannot expect to be able to compute a border basis
for \emph{any} order ideal, simply as such a basis
does not necessarily exist.
As \emph{a priori} it is unclear which are the admissible order ideals
\(\oi\)
for an ideal \(I\) given by generators, we use an indirect way to
specify \(\oi\):
we use a weight vector \(w \in \R^{\bT^{n}}\)
and want to find \(\oi \in \Lambda(I)\)
maximizing the total weight \(\sum_{m \in \oi} w_{m}\) of
\(\oi\). Note that any admissible order ideal \(\oi\) can be specified
via an appropriate weight vector \(w\) so our approach, while
indirect, is without loss of generality. 
As \(w\) is an infinite vector, in practice it should be probably
given explicitly for a finite number of coordinates, and the remaining
coordinates are declared to be \(0\) or some other fixed value; this
is not a restriction as all admissible order ideals are finite and the
occuring maximum degree is bounded. 
Recall that \(\Lambda(I)\) denotes the set of all admissible order
ideals of \(I\).
We will show
how to compute such a weight-maximal $\oi$ and its border basis
for a zero-dimensional ideal $I \subseteq \pring$.

We adapt the border basis algorithm
in \cite{kehrein2006cbb}.
\begin{algorithm}{Generalized border basis algorithm}{\bbasis}%
{$F$ a finite generating set of a zero-dimensional ideal,
and a weight vector \(w\) on \(\bT^{n}\).}%
{$\mathcal{G}$ a border basis of the ideal.}
\label{alg:borderBasisNew}

\item Let $d\coloneqq \max_{f\in F} \deg(f)$.
\item \label{enu:bbasisStep2_1}
  $M \coloneqq \lstabspan(F, \bT^{n}_{\leq d})$
  using a degree-compatible ordering on \(\bT^{n}_{\leq d}\)
  \\ 
  (i.e., \(m_{1} < m_{2}\) whenever \(\deg m_{1} < \deg m_{2}\)).
\item \label{enu:bbasisStopNew}
If \(\bT^{n}_{= d} \nsubseteq \LT(M)\)
then set $d\coloneqq d+1$
and
go to step \ref{enu:bbasisStep2_1}.
\item \label{enu:bbasisDimension}
  Set \(d_{\text{old}} \coloneqq d\),
  \(d \coloneqq \size{\bT^{n}_{\leq d}} - \lvert M \rvert\).
  If \(d \leq d_{\text{old}}\) then let \(M \coloneqq M^{\leq d}\).
  Otherwise let \(M \coloneqq \lstabspan(F,\bT^{n}_{\leq d})\).
  \item \label{enu:calcOrderIdeal}
    Write up the system \eqref{eq:OrderIdealPolytope}
    of inequalities for \(M\) and \(d\).
    Choose an integral solution \(z\) maximizing \(w z\).
    Set $\oi \coloneqq \{ m \in \bT^{n}_{\leq d - 1} : z_{m} = 1\}$.
  \item \label{enu:basisTransformation} Let $\mathcal{G} \coloneqq   \basisT(M,\oi)$.
  \end{algorithm}
Our \emph{generalized border basis algorithm} \ref{alg:borderBasisNew}
first determines the right computational
universe $\bT^n_{\leq d}$ until step \ref{enu:bbasisDimension},
i.e., a large enough $d \in \N$
such that the associated
$\bT^{n}_{\leq d}$-stabilized span $M$ contains all border bases.
Here step \ref{enu:bbasisStopNew}
is a convenient way to quickly check whether
the universe is already large enough.
Step \ref{enu:bbasisDimension} adjusts \(d\) to
the actual dimension of \(\pring / I\)
and adjusts \(M\).

In the second phase,
step~\ref{enu:calcOrderIdeal} optimizes over
the order ideal polytope $P(\igen{M})$
to find an optimal admissible order ideal
using a mixed integer programming solver,
and then step~\ref{enu:basisTransformation}
computes the corresponding border basis.
The main idea of this last step is to
apply Gaussian elimination to $M$
to bring it into a form where with the exception of the leading terms,
all monomials are from \(\oi\).

\begin{lem}
  \label{lem:basisT}
  Let $L = \bT^n_{\leq \ell}$ with $\ell\in \N$,
  let $M$ be a non-empty finite set of polynomials satisfying
  \(\vgen{M} = \igen{M} \cap \vgen{L}\)
  and
  let $\oi$ be an order ideal with
  $\partial \oi \subseteq L$
  and $\oi \in \Lambda(\igen{M})$.
  Then Algorithm~\ref{alg:basisT} returns
  an $\oi$-border basis $\mathcal G$ of $\igen{M}$.
\begin{proof}
First, because \(M \subseteq L\) and \(M\) is non-empty,
clearly the largest degree is \(\ell\) among the polynomials in
\(\vgen{M} = \igen{M} \cap \vgen{L}\) and hence in \(M\).
Thus step \ref{enu:l} computes the correct value of \(\ell\).

As $\oi \in \Lambda(\igen{M})$
we have $\pring = \igen{M} \oplus \vgen{\oi}$
and hence
\begin{equation*}
  \vgen{L}
  =
  \vgen{L} \cap \left( \igen{M} \oplus \vgen{\oi} \right)
  =
  \left( \vgen{L} \cap \igen{M} \right) \oplus \vgen{\oi}
  =
  \vgen{M} \oplus \vgen{\oi}
\end{equation*}
by the modular law,
in particular,
\(\size{\mathcal{G}'} = \size{M} = \size{L} - \size{\oi}\).
Now, none of the polynomials in \(\mathcal{G}'\) are supported
on \(\oi\), and as \(\oi\) is an initial segment of the ordering
used for Gaussian elimination,
it follows that all the leading terms in \(\mathcal{G}'\)
lie in \(L \setminus \oi\).
Since \(\size{\mathcal{G}'} = \size{M} = \size{L} - \size{\oi}\),
it follows that all \(m \in L \setminus \oi\) appear as leading term
exactly once in \(\mathcal{G}'\), and hence not as other term,
i.e., all polynomials \(g \in \mathcal{G}'\) have the form
\begin{equation*}
  g = m_{0} - \sum_{m \in \oi} a_{m} m, \qquad m_{0} \in L \setminus \oi.
\end{equation*}
Obviously, restricting to the polynomials where the
leading term is a border element of \(\oi\)
in step~\ref{enu:restrictToBorder}
provides a border basis of \(\oi\).
\end{proof}
\end{lem}

\begin{algorithm}{Basis transformation algorithm}{\basisT}%
{$M,\oi$ as in Lemma \ref{lem:basisT}.}%
{$\mathcal G$ as in Lemma \ref{lem:basisT}.}
\label{alg:basisT}

\item\label{enu:l}
  Set \(\ell \coloneqq \max_{m \in M} \deg(m)\).
\item
  Reduce \(M\) using Gaussian elimination
  (Algorithm~\ref{alg:gaussEl})
  using an ordering where \(\oi\) is an initial segment
  (i.e., consists of the smallest elements):
  \(\mathcal{G}' \coloneqq \gausEl(M)\).
\item\label{enu:restrictToBorder}
  Return $\mathcal{G} \coloneqq
  \{g \in \mathcal G' : \LT(g) \in \partial\oi \}$.
\end{algorithm}

We will show now that Algorithm~\ref{alg:borderBasisNew} computes an $\oi$-border basis for $\oi \in \Lambda(I)$.

\begin{prop}
\label{pro:bbasisNew}Let $F \subseteq \pring$
be a finite set of polynomials that generates a zero-dimensional ideal
$I=\igen{F}$. Then Algorithm \ref{alg:borderBasisNew} computes
the $\oi$-border basis $\mathcal{G}$ of $I$ for any (chosen) $\oi \in \Lambda(I)$.
\begin{proof}
Till step~\ref{enu:bbasisStopNew},
the algorithm step by step enlarges the computational universe
\(\bT^{n}_{\leq d}\) via increasing \(d\).
Since \(I\) is zero-dimensional,
the test \(\bT^{n}_{=d} \subseteq \LT(M)\) will be true for large
enough \(d\), hence the algorithm will eventually reach
step~\ref{enu:bbasisDimension}.
By Corollary~\ref{cor:idealApprox}, we have then
\(\vgen{M} = I^{\leq d}\),
and that
step \ref{enu:bbasisDimension} sets \(d\) to
the dimension of \(\pring/I\).
It also updates \(M\) so that together with the new \(d\)
it satisfies \(\vgen{M} = I^{\leq d}\).
Obviously, \(\bT^{n}_{\leq d}\) contains all
order ideals supporting a border basis, i.e., all $\oi \in \Lambda(I)$
and even the boundary of these order ideals.
Observe that $I = \igen{F} = \igen{M}$ and thus, by
Lemma~\ref{lem:basisT}, it follows that $\mathcal G$ is indeed an
$\oi$-border basis of $\igen{F}$.
Note that
$\partial \oi \subseteq \bT^{n}_{\leq d}$
follows
via \(\left. \vgen{M}^{\leq d} / \vgen{M}^{\leq d-1} \right. \cong
 \vgen{\bT^{n}_{= d}}\)
as \(d\) is the dimension of \(\pring / I\).

We conclude that \(M\) satisfies the conditions of
Lemma~\ref{lem:polytope-comp}, e.g.,
\(d = \size{\bT^{n}_{\leq d}} - \size{M}\)
is ensured by step~\ref{enu:bbasisDimension},
and therefore the integral solutions of the system used
in step~\ref{enu:calcOrderIdeal} are all the characteristic
vectors of admissible order ideals.
This step therefore really computes
a maximal-weight admissible order ideal \(\oi\).
Finally, step~\ref{enu:basisTransformation}
returns a border basis of \(\oi\) by Lemma~\ref{lem:basisT}.
\end{proof}
\end{prop}

The border basis algorithm in \cite{kehrein2006cbb}
allows using computational universes $L$
smaller than \(\bT^{n}_{\leq d}\),
improving performance of the algorithm.
However, as we want to consider \emph{all}
order ideals and border bases,
we deliberately chose the computational universe large enough
to contain all possible order ideals.
If a subset of all admissible
order ideals is sufficient,
then the same optimizations can be applied throughout.

For certain choices of the weight vector $w$ though
it can be hard to compute a
maximum weight order ideal as we will show now. In fact this also
shows that there is no general, efficient way of specifying any admissible order
ideal (unless \(\coNP = \NP\)).

\section{Complexity of finding maximum weight order ideals}
\label{sec:complexity-finding-ord-ideal}

In this section,
we show that finding a maximum weight, admissible order ideal of
a zero-dimensional ideal given by generators is \NP-hard
(Theorem~\ref{thm:MaxOrderIdeal-NPhard}).
The hardness result is
unexpected in the sense that we merely ask for a \emph{nice} basis transformation. On the
other hand it highlights the crucial role of order ideals in describing the
combinatorial structure of the ideal. As an immediate consequence it follows that it is rather unlikely that we can
obtain a good characterization of the integral hull
of the order ideal polytope $P(\igen{M})$ (unless \(\NP = \coNP\))
and we will not be able to compute order ideals that support
a border basis and have maximum weight efficiently in
the worst case (unless $\NP = \Pclass$). This shows
that it is hard
not only to compute the necessary liftings of the initial set of
polynomials via the $\lstabspan$ procedure but also to actually
determine an optimal choice of an order ideal once an $L$-stable span has
been computed.

From a practical point of view this is not too problematic as,
although \NP-hard, computing a maximum weight order ideal is no harder
than actually computing the $\lstabspan$ in general. For bounds on the
degree $d \in \N$ needed to compute border bases, see e.g.,
\cite[Lemma 2.4]{deLoera2009ecps}; the border basis algorithm
generates the Nullstellensatz certificates and is therefore subject to
the same bounds. Further, state-of-the-art mixed integer programming
solvers such as \textprogram{scip} \cite{achterberg2009scip},
\textprogram{cplex} \cite{cplex200811}, or \textprogram{gurobi}
\cite{gurobi} can handle instance sizes far beyond the point for which
the actual border bases can be computed. Very good solutions can also
be generated using simple local search schemes starting from a
feasible order ideal derived from a degree-compatible term ordering.

\subsection{Fast without constraint}
\label{sec:fast-with-constr}

\begin{figure}
  \centering
\includegraphics[scale=0.55]{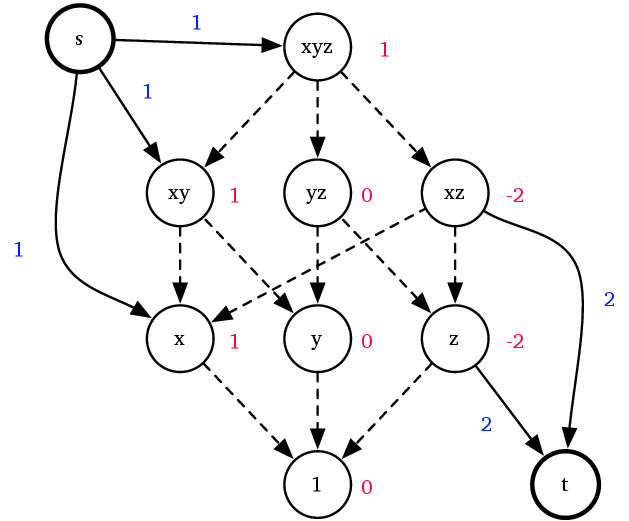}
\caption[Order ideal computation as minimum cut problem]%
{\label{fig:orderIdeal}%
Order ideal computation as minimum cut problem.
In this example, the order ideal consists of all monomials dividing
\(xyz\), i.e., \(\{1, x, y, z, xy, yz, xz, xyz\}\),
with (arbitrarily) chosen weights \(\{0, 1, 0, -2, 1, 0, -2, 1\}\).
The weight \(w_{u}\) of a monomial \(u\) is shown in red
next to the node of the monomial.
Arcs are labelled in blue with their capacity resulting
from the weights of monomials.
Arcs with capacity $\infty$ are dashed.
Arcs with capacity \(0\) are omitted.}
\end{figure}

Determining an order ideal of maximum
weight (not necessarily supporting a border basis!) in a computational universe \(L\) without having any
constraints on the dimension of the respective spaces can be done in time
polynomial in $\lvert L \rvert $ as we will show now. This follows
with \cite{picard1976maximal} and we simply transform the maximum
weight order ideal problem into a \emph{minimum cut problem}. For this let $w \in \Z^L$ be a
weight vector. We define a
directed
graph $\Gamma \coloneqq (V,A)$ with $V \coloneqq L \cup \{s,t\}$ and $\tilde A \coloneqq  \{ (u,v)
\mid u,v \in L\text{ and } v \mid u \}$, i.e., whenever $v \mid u$ we add an arc from
$u$ to $v$.
In fact, it is enough to have an arc when \(u = v x\) for some variable \(x\), i.e., to consider the transitive reduction of $\tilde A$.
Define \[A \coloneqq  \tilde A \cup \{ (s,u) \mid u \in L, w_u > 0\} \cup
\{ (u,t) \mid u \in L, w_u < 0 \}.\]
Now we turn to the arc capacities.
Let \(\kappa(u,v)\) denote the capacity or arc \((u,v)\)
defined as follows.
For \(u\) and \(v\) both in \(L\),
we set \(\kappa(u,v) \coloneqq \infty\).
We set \(\kappa(s,v) \coloneqq w_{v}\)
and \(\kappa (u,t) = - w_{u}\)
for \(u\), \(v\) in \(L\).
An example is depicted in \autoref{fig:orderIdeal}.

For $U,W \subseteq V$, we define $C(U,W) \coloneqq \sum_{(u,w) \in U
  \times W} \kappa(u,w)$ as the directed cut value.  An \emph{(s,t)-cut} $(S,\bar S)$ is a
partition $S \stackrel{\cdot}{\cup} \bar S=V$ of the vertices of $V$
with $s \in S$ and $t \in \bar S$ and the \emph{weight} of the cut is
$C(S,\bar S)$; note that the direction of the arcs matters. We would like to compute an order ideal contained in $L$ with maximum weight:

\[ \max \left\{ \sum_{u \in \oi} w_u \middle| \oi \subseteq L
  \text{ order ideal} \right\}.\]
Observe that $(S, \bar{S})$ is a directed cut in $\Gamma$ of finite weight,
if and only if
there exists no arc $(u,v) \in \tilde A$ with $u \in S$ and
$v \in \bar S$, i.e., for all monomials \(u, v \in \oi\)
with \(v \mid u\), if \(u \in S\) then \(v \in S\).
In other words,
$(S, \bar{S})$ is a cut in $\Gamma$ of
finite weight if and only if $S \setminus \{s\}$ is an order ideal. We
can therefore rewrite the optimization problem as follows:
\begin{align*}
  \MoveEqLeft
    \max \left\{
      \sum_{u \in \oi} w_u \middle| \oi \subseteq L \text{ order
        ideal} \right\} \\
    &
    = \max \{ C(\{s\},\oi ) - C(\oi, \{t\}) \mid
    \oi \subseteq L \text{ order ideal}\} \\
    &
    = \max \{ C(\{s\}, L ) - C(\{s\},L \setminus \oi ) - C(\oi, \{t\}) \mid
    \oi \subseteq L \text{ order ideal} \} \\
    &
    = C(\{s\}, L ) - \min \{ C(\{s\},L \setminus \oi ) + C(\oi, \{t\}) \mid
    \oi \subseteq L \text{ order ideal} \}  \\
    &
    = C(\{s\}, L ) - \min \{ C(\{s\} \cup \oi , (L \setminus \oi) \cup \{t\}
    ) \mid \oi \subseteq L \}.
\end{align*}
The last line asks for a minimum weight cut in the graph
$\Gamma$. Note that we can indeed drop the condition that $\oi$ has to
be an order ideal as it is guaranteed implicitly by all finite weight
cuts as explained above.  The minimum cut can now be computed in
polynomial time in the number of vertices and arcs (see e.g.,
\cite{wolsey1999integer}) and so can an order ideal $\oi$ of maximum
weight efficiently.

\subsection{\NP-hard with constraints}
\label{sec:np-hard-constraints}

So far we did not include the additional requirements as specified by
the order ideal polytope (see \eqref{eq:OrderIdealPolytope}), in
order to obtain order ideals that do actually support a border basis
of the ideal $I$ under consideration. We will now show that when
including these additional requirements, the problem of computing an
order ideal of maximum weight becomes \NP-hard.
In \cite[Discussion after
Definition 3.2]{hochbaum2000performance} it was indicated that
determining a maximum weight order ideal of a pre-defined size is
\NP-hard by a reduction from MaxClique, however this is
different from our problem, as we have additional constraints coming
from the dimension of the factor spaces of the ideal (see
constraints~\eqref{eq:equiSteinitz}).

We will show \NP-hardness by a reduction from the
\problemref{pr:kClique} problem, which is well known to be \NP-complete
(see, e.g., \cite{garey1979cai} or
\cite[GT22]{crescenzi-compendium}).
Given an undirected simple graph $\Gamma = (V,E)$,
recall that a \emph{clique} $C$ is a subset of $V$ such that for all
distinct $u,v \in C$ we have $(u,v) \in E$. We consider the decision problem: 

\begin{problem}[$k$-Clique]
  \label{pr:kClique}
  Let $\Gamma = (V,E)$ be an undirected simple graph.
  Decide whether $\Gamma$ contains a
  clique $C$ of size $k$.
\end{problem}

Our optimization problem of interest is:

\begin{problem}[Maximum weight admissible order ideal]
  \label{pr:MaxOrderIdeal}
  Let $M \subseteq \pring$ be a system of polynomials generating a
  zero-dimensional ideal and let $w\in \Z^{\bT^n}$ be a weight
  on the monomials.
  Compute
  an admissible order ideal $\oi \subseteq \bT^{n}$ for $\igen{M}$
  with maximum weight \(\sum_{m \in \oi} w_{m}\)
  with respect to $w$, i.e., compute
  \begin{equation*}
    \argmax_{\oi \in \Lambda(\igen{M})}
    \sum_{m \in \oi} w_{m}.
  \end{equation*}
\end{problem}

By a reduction from \problemref{pr:kClique} we obtain:

\begin{thm}
\label{thm:MaxOrderIdeal-NPhard}
\problemref{pr:MaxOrderIdeal} is \NP-hard over ground fields \(K\)
of characteristic \(0\).
\end{thm}

As a preparation for the proof, we show that
for every graph $\Gamma = (V,E)$ and $k\in [\lvert V
\rvert]$ there exists a
system of polynomials
$F_{\lvert V \rvert,k} \subseteq K[x_v \mid v \in V]$
spanning a zero-dimensional ideal
such that solving the \problemref{pr:MaxOrderIdeal} problem for
$F_{\lvert V \rvert,k}$ solves the \problemref{pr:kClique} problem for $\Gamma$.
For this, we construct an ideal encoding all $k$-cliques of the
complete graph on $n$ vertices: Let $n \in \N$ and $k \in [n]$ and define 
\begin{equation*}
  F_{n,k} \coloneqq \{ v_j \mid j \in [n-k] \} \cup \bT^n_{=3}
\end{equation*}
with $ v_j \coloneqq \sum_{i \in [n]} i^j x_i$.
We consider the ideal generated by \(F_{n,k}\).
We show that its order ideals are in one-to-one correspondence with the
\(k\)-element subsets of the set of \(n\) variables \(x_1,\dots,x_n\)
as stated in the following
lemma.

\begin{lem}
\label{thm:thePolySystem}
Let \(K\) be a field of characteristic \(0\) together with
$n \in \N$ and $k\in [n]$.
Then $F_{n,k}$ generates a zero-dimensional
ideal such that $\oi \in \Lambda \left( \igen{F_{n,k}} \right)$ if and only if
$\oi^{=1} \subseteq \bT^n_{=1}$ with $\lvert \oi^{=1} \rvert = k$, $\oi^{=2} = \{xy
\mid x,y \in \oi^{=1}\}$, and $\oi^{=\ell} = \emptyset$ for all $\ell \geq 3$.
\begin{proof}
We start by providing an explicit representation of the factor ring
$\left. K[x_{1}, \dots, x_{n}]
  \middle/ \igen{F_{n,k}} \right.$.
As \(F_{n,k}\) consists of homogeneous polynomials,
it generates a homogeneous ideal \(I\),
and induces a degree decomposition of the factor ring:
\begin{equation*}
  \frac{K[x_{1}, \dots, x_{n}]}{\igen{F_{n,k}}}
  = \bigoplus_{i=0}^{\infty} \frac{\vgen{\bT^{n}_{=i}}}{I^{=i}}.
\end{equation*}
To actually determine the factors,
let \(x_{i_{1}}\), \dots, \(x_{i_{k}}\) be \(k\) many
distinct variables from \(x_{1}\), \dots, \(x_{n}\).
We prove that
\(x_{i_{1}}\), \dots, \(x_{i_{k}}\), \(v_{1}\), \dots, \(v_{n-k}\)
is a vector space basis of \(\vgen{\bT^{n}_{=1}}\) by showing
that its coefficient matrix
in the standard basis \(x_{1}\). \dots, \(x_{n}\)
has non-zero determinant.
Expanding the determinant by the \(k\) rows of
\(x_{i_{1}}\), \dots, \(x_{i_{k}}\),
each of which contains only one non-zero element,
the determinant becomes equal to up to a sign
to the Vandermonde matrix of the set of numbers
\([n] \setminus \{i_{1}, \dots, i_{k}\}\),
and hence it is indeed non-zero.

The ring \(\pring\) is also a polynomial ring
in any basis of \(\vgen{\bT^{n}_{=1}}\),
and the basis \(x_{i_{1}}\), \dots, \(x_{i_{k}}\), \(v_{1}\),
\dots, \(v_{n-k}\) is particularly suitable to determine
the factor \(\pring / I\) together with the degree decomposition:
\begin{equation}
  \label{eq:Vandermonde-decomposition}
 \begin{split}
  \frac{K[x_{1}, \dots, x_{n}]}{\igen{F_{n,k}}}
  &
  =
  \frac{K[x_{i_{1}}, \dots, x_{i_{k}}, v_{1}, \dots, v_{n-k}]}
  {\igen{v_{1}, \dots, v_{n-k}, \bT^{n}_{=3}}}
  \\
  &
  =
  \frac{K[x_{i_{1}}, \dots, x_{i_{k}}]}
  {\igen{\bT^{n}_{=3}}}
  =
  \vgen{1} \oplus \vgen{x_{i_{1}}, \dots, x_{i_{k}}}
  \oplus
  \vgen{x_{i} x_{j} : i,j \in \{i_{1}, \dotsc, i_{k}\}}
 \end{split}
\end{equation}
where the generating sets are actually bases
of the respective degree components.

Given an order ideal \(\oi\) of \(I\),
the isomorphism \(\vgen{\oi} \cong \pring / I\)
clearly preserves the degree decomposition,
i.e.,
\begin{math}
  \vgen{\oi^{=\ell}} \cong ( \pring / I )^{=\ell}
\end{math}
for all \(\ell\).
Hence \(\size{\oi^{=0}} = 1\),
\(\size{\oi^{=1}} = k\), \(\size{\oi^{=2}} = \binom{k+1}{2}\),
and \(\size{\oi^{=\ell}} = 0\) for \(\ell \geq 3\).
It follows that \(\oi\) has the claimed form,
in particular,
\(\oi^{=2} = \{xy \mid x, y \in \oi^{=1}\}\)
as the left-hand side is clearly a subset of the right-hand side,
and they have the same finite size.

For the other direction,
let \(\oi^{=1} = \{x_{i_{1}}, \dots, x_{i_{k}}\}\) with
\(\oi^{=2} = \{xy \mid x, y \in \oi^{=1}\}\),
\(\oi^{=0} = \{1\}\)
and \(\oi^{=\ell} = \emptyset\) for \(\ell \geq 3\).
Then \(\oi\) is an order ideal,
and \eqref{eq:Vandermonde-decomposition} shows
that the image of \(\oi\) in \(\pring / I\) is a basis.
Thus \(\oi\) is an admissible order ideal for \(I\), as claimed.
\end{proof}
\end{lem}

Note that the order ideals of $F_{n,k}$ indeed correspond to the
$k$-cliques of the complete graph on $n$ vertices: If $\oi \in
\Lambda(F_{n,k})$, then $\oi^{=1} = \{x_{i_1}, \dots, x_{i_k}\}$ and
$x_{i_j}x_{i_l} \in \oi^{=2}$ if and only if $x_{i_j},x_{i_l} \in \oi^{=1}$. If we now
remove all elements of the form $x_{i_j}^2$ with $x_{i_j} \in \oi^{=1}$, and
there are $k$ of those, then \[ \lvert \oi^{=2} \setminus \{x_{i_j}^2 \mid
x_{i_j} \in \oi^{=1}\} \rvert = \frac{k(k-1)}{2},\]
the size of a $k$-clique. We are ready to prove the main result of
\hyperref[sec:Computing-border-bases]{this section}.

\begin{proof}[Proof of Theorem~\ref{thm:MaxOrderIdeal-NPhard}]
The proof is by a reduction from the \NP-hard \problemref{pr:kClique}
problem.
Let us start with an instance of \problemref{pr:kClique},
i.e., an undirected graph $\Gamma = (V,E)$
with $n \coloneqq \lvert V \rvert$ and $k \in [n]$.
We consider $M \coloneqq F_{n,k}$ and
define $w \in \Z^{\bT^n}$ via
\[w_m =
\begin{cases}
1, &\text{if } m= x_u x_v \text{ and either } (u,v) \in E \text{ or } u=v; \\
0, &\text{otherwise},
\end{cases}
\]
for all \(m \in \bT^{n}\).
By Lemma~\ref{thm:thePolySystem},
there is a bijection of the admissible order
ideals \(\oi\) of $\igen{M}$
and the $k$-cliques of the complete graph on $n$ vertices
given by
\[C_{\oi} \coloneqq \{v \in V \mid x_v \in \oi \}.\]
The weight of \(\oi\) is the sum of
the weights of the monomials \(x_{u} x_{v}\) in \(\oi\).
To the weight of \(\oi\),
the contribution of the monomials with \(u =
v\), i.e., of the form \(x_{v}^{2}\) is
the number of vertices of \(C_{\oi}\), i.e., \(k\).
The monomials \(x_{u} x_{v}\) with \(u \neq v\)
contribute the number of edges in \(C_{\oi} \cap \Gamma\)
to the weight of \(\oi\).
Hence the weight of \(\oi\) is the sum of \(k\)
and the number of edges in \(C_{\oi} \cap \Gamma\).

The largest possible value of this weight is $k(k+1)/2$,
and this is realized exactly by cliques \(C_{\oi}\) of \(\Gamma\)
of size \(k\).
(If such cliques do not exist, then the maximal weight
is less than \(k (k+1) / 2\).)
All in all, the maximum weight is \(k (k+1) / 2\)
if and only if $\Gamma$
contains a clique \(C_{\oi}\) of size $k$.
We obtain that
\problemref{pr:MaxOrderIdeal} solves \problemref{pr:kClique} and so
the former has to be \NP-hard.
\end{proof}

\subsection{Extension complexity of admissible order ideals}
\label{sec:extens-compl-order}

The order ideal polytope \(P(I)\) was introduced as a relaxation
of the convex hull \(\OIP(I)\) of (the characteristic vectors of)
all admissible order ideals of the ideal \(I\).
Therefore one might wonder whether there exists a description with a
polynomial number of
linear inequalities 
of the convex hull \(\OIP(I)\).
This question is the natural counterpart of algorithmic complexity
in the context of linear programming.
Here we show that in general \(\OIP(I)\) requires a subexponential
number of inequalities in the size of the computational universe,
even if one allows additional extra variables,
i.e., the \emph{extension complexity} (see below) of \(\OIP(I)\) is
subexponential. As customary in extended formulations this result does
not depend on any complexity theoretic assumptions,
see \cite{ConfortiCornuejolsZambelli10,Kaibel11,extform4,BPZ2014}
for details.
The result could be also formulated
independent of the order ideal polytope,
namely, that the linear programming formulation complexity
(complexity measured in the size of a linear program)
of the combinatorial problem
to find a maximum-weight admissible order ideal is
subexponential.
However, for simplicity, we stick to the polyhedral formulation,
and refer the
interested reader to
\cite{BPZ2014} for the general model.

Recall that the \emph{extension complexity} \(\xc(P)\)
of a polyhedron \(P\) is the minimum number of facets
of a polyhedron \(Q\), such that \(P\) is an affine image of
\(Q\). The extension complexity captures the inherent complexity of a
polytope being expressed by means of linear inequalities.

\begin{thm}
  \label{thm:xc-admissible}
  For any ground field \(K\) of characteristic \(0\),
  there is an ideal \(I\) of \(K[x_{1}, \dots, x_{2n}]\)
  such that all admissible order ideals of \(I\)
  contain monomials only up to degree \(2\), and
  \begin{equation}
    \xc(\OIP(I)) = 2^{\Omega(n)}.
  \end{equation}
\begin{proof}
We shall use the ideal \(I\) generated by \(F_{2n,n}\)
from Lemma~\ref{thm:thePolySystem},
whose admissible order ideals have a nice description,
and all of which consist of monomials only up to degree \(2\).
Therefore we obtain the following description of \(\OIP(I)\):
\begin{align}
  \OIP(I) &= \conv{\{y^{S} : S \subseteq [2n], \size{S} = n\}}
  \subseteq [0,1]^{\bT^{n}_{=2}}
  \\
  y^{S}_{x_{i} x_{j}} &=
  \begin{cases}
    1, & \text{if } i,j \in S, \\
    0, & \text{otherwise}.
  \end{cases}
\end{align}
Here for simplicity we restrict
to the relevant coordinates only, the other coordinates are
affine combinations of these (e.g., \(y_{x_{i}} = y_{x_{i}^{2}}\),
\(y_{1} = 0\)).
We will show that the correlation polytope
is an affine projection of \(\OIP(I)\),
and therefore \(\xc(\OIP(I)) \geq \xc(\CORR(n)) = 2^{\Omega(n)}\)
by \cite[Lemma~9(i)]{extform4} and
\cite[Theorem~4(i)]{bfps2012jour}.

Recall that the \emph{correlation polytope}
is the convex hull of all 0/1-matrices of rank \(1\):
\begin{align}
  \CORR(n) &\coloneqq
  \conv{\{\allOne_{S} \allOne_{S}^{\intercal} \mid S \subseteq [n]\}}
  \subseteq [0,1]^{n \times n},
  \\
  \allOne_{S}(i) &\coloneqq
  \begin{cases}
    1, & \text{if } i \in S, \\
    0, & \text{otherwise}.
  \end{cases}
\end{align}
An affine projection \(f \colon \OIP(I) \to \CORR (n)\)
is clearly provided by
\begin{align}
  f(y)_{i,j} &\coloneqq y_{x_{i} x_{j}} & i,j &\in [n], \\
  \intertext{where the vertices of \(\OIP(I)\) are mapped to
    vertices of \(\CORR(n)\)}
  f(y^{S}) &= \allOne_{S \cap [n]} \allOne_{S \cap [n]}^{\intercal}
  & S &\subseteq [2n], \size{S} = n.
\end{align}
Note that \(2n\)
variables were chosen for \(\OIP(I)\)
so that every subset \(T \subseteq [n]\) arises as an intersection
\(T = S \cap [n]\) for some \(S \subseteq [2n]\) of size \([n]\).
\end{proof}
\end{thm}

\subsection{Discussion of the complexity of
  finding maximal weight admissible order ideals}
\label{sec:compl-order-ideal}

We now briefly summarize the implications of these
complexity results.
Note that the hardness proof in
\autoref{sec:np-hard-constraints} is \emph{independent} of the
order ideal polytope,
and shows worst-case hardness for any algorithm.

\begin{enumerate}[itemsep=1em]
\item \emph{No general characterization of all admissible order
    ideals.} The hardness in
  \autoref{sec:np-hard-constraints} is established for the
  \problemref{pr:MaxOrderIdeal} problem.
  As such, unless \(\NP = \coNP\), which
  is generally believed to be not the case, in general there will be no
  \emph{good characterization} of order ideals that will be admissible
  for a given ideal. Complementing this, the result in
  Section~\ref{sec:extens-compl-order} rules out any small linear
  programming formulation for the convex hull of
  admissible order ideals
  irrespective of \(\NP\) vs.~\(\coNP\). 

\item \emph{No theoretically efficient algorithm for computing
    maximal weight order ideals.}
  In particular, unless \(\NP = \Pclass\), there
  will be no polynomial time algorithm computing a maximal weight
  order ideal.  However, this is worst-case complexity,
  and does not necessarily capture well real-world
  performance,
  as e.g., the Traveling Salesman Problem is also NP-hard,
  however solvable for
  real-world instances with millions of cities
  in reasonable computational time
  (see e.g., \cite{applegate2011traveling}).

\item \emph{Real-world computational complexity.}  While the
  determinination of a maximal weight admissible order ideal is
  theoretically NP-hard as discussed above, in practice this problem
  can be solved very easily with state-of-the-art solvers
  such as e.g.,
  \texttt{scip}, \texttt{CPLEX}, or \texttt{Gurobi}, typically in the
  order of seconds. Comparing the generalized border basis algorithm
  (Algorithm~\ref{alg:borderBasisNew}) to the border basis algorithm in
  \cite{kehrein2006cbb},
  the major difference is the additional computational steps
  \ref{enu:calcOrderIdeal} and \ref{enu:basisTransformation}.
  The basis transformation in \ref{enu:basisTransformation} is
  very cheap, and so is 
  step \ref{enu:calcOrderIdeal} for all practical purposes as
  indicated.

The real bottleneck in our border basis algorithm
  (and also the one in \cite{kehrein2006cbb}, upon which ours is
  based) is computing the \(L\)-stable
  span, which can be several orders of magnitude slower than
  determining the order ideal. In summary, we believe that our method
  has little additional costs compared to the border basis algorithm
  in \cite{kehrein2006cbb}, when incorporated correctly into a
  state-of-the-art implementation.

\item \emph{Size of the order ideal polytope.}
  The description of the
  order ideal polytope in Definition~\ref{def:orderIdealPolytopeAlg}
  has a number of inequalities of roughly \(O(n^{d^{2}})\) due to
  \eqref{eq:steinitz2}.
  This is roughly the largest possible number of facets
  of the order ideal polytope,
  i.e., the number of all order ideals.
  Recall that the order ideal polytope is not necessarily
  the convex hull of all admissible order ideals, but only a
  relaxation of it, and as Theorem~\ref{thm:xc-admissible} shows
  (where \(d=2\)), the convex hull itself
  requires much more inequalities in the worst case.
  However, it is
  conceivable that the convex hull admits an even smaller relaxation
  than the order ideal polytope
  via uncapacitated network flows
  or separating the inequalities
  (see e.g., \cite{schrijver1986theory})
  as e.g., done for the spanning tree polytope.
\end{enumerate}

\section{\label{sec:Computational-results}Computational results}

We performed computational tests to verify the practical feasibility
of our method,
with an emphasis of the optimization step
over the order ideal polytope, once it is written down,
as this is the new aspect in our algorithm.
Because this is not expected to be the bottleneck,
we refrained from a comprehensive performance test,
and used small-sized problems.

For simplicity, we computed only degree-compatible
order ideals.  All computations were performed with \textprogram{CoCoA
  4.7.5} \cite{team2009cocoa} and \textprogram{scip 1.1.0}
\cite{achterberg2009scip} on a 2 Ghz Dual Core Intel machine with 2
GB of main memory\footnote{Source code is available at:
  \url{https://app.box.com/s/fcxtocvpqqj0b2dezj40v4tn24sfkk1m}}. 

\subsection*{Test setup}
\label{sec:test-setup}

The employed methodology was as
follows. We first computed a border basis using the border
basis algorithm in \cite{kehrein2006cbb}.
From the last run of the algorithm we extracted the
$L$-stabilized span and brought it into canonical form as the actual
\(L\)-stable span computation is not the focus here but the
computation of admissible order ideals. We then generated
the constraint \eqref{eq:degreeCompatibility} from the order ideal
that we obtained; from the $L$-stabilized span in matrix from, we
generated the constraints \eqref{eq:equiSteinitz} adapted to
degree-compatible order ideals.
We performed computations on various sets of
systems of polynomial equations. We then transcribed these
constraints into the CPLEX LP format which served as input for
\textprogram{scip}. For the optimization we chose various weight
vectors.  We tested random weight vectors and we constructed a
weight vector with the intent to make the optimization particularly
hard by giving monomials deep in the order ideal negative weights and
assigning positive weights for the outer elements.

\subsection*{Results}
\label{sec:results}

We report the results of our tests in Table~\ref{tab:results}. In all
cases, the optimization (i.e., the computation of the maximum weight
order ideal) was performed in less than a second
(see column \emph{optimization}), whereas the actual
calculation of the initial border bases was significantly more time
consuming. This is not unexpected as
the computation of the $L$-stable span is significantly more involved
than computing a maximum weight order ideal: the former can be double
exponential whereas the latter is at most single exponential via
complete enumeration.

\subsection*{An example application: counting order ideals}
\label{sec:an-example-appl}

When computationally feasible,
we also counted all feasible order ideals with
\textprogram{scip}, which basically means enumerating all feasible
solutions, to demonstrate feasibility for reasonably sized
instances. This is reported in column \emph{counting}.

\begin{table}[htp]
\begin{center}
\begin{tabular}{|>{\centering}p{4cm}|c|c|c|c|}
\hline
 polynomial system &
 \begin{tabular}{c}
   degree vector of\\ order ideal
 \end{tabular}
 &
 optimization [s] & counting [s] & \# order ideals \\
\hhline{|=|=|=|=|=|}
$x^3, xy^2+y^3$ & $(1,2,3,2,1)$ & < 0.01 & 0.02 & 3 \\
\hline
vanishing ideal of the points \((0,0,0,1)\), \((1,0,0,2)\), \((3,0,0,2)\),
\((5,0,0,3)\), \((-1,0,0,4)\), \((4,4,4,5)\), \((0,0,7,6))\).
 & $(1,4,2) $ & < 0.01 & 0.02 & 45 \\
\hline
$x+y+z-u-v$, $x^2-x$, $y^2-y$, $z^2-z$, $u^2-u$, $v^2-v$ & $(1,4,5) $ & < 0.01 & 0.35 & 1,260 \\
\hline
$x+y+z-u-v$, $x^3-x$, $y^3-y$, $z^2-z$, $u^2-u$, $v^2-v$ & $(1,4,7,6)$ & 0.02 & 51.50 & 106,820 \\ 
\hline
$x+y+z-u-v$, $x^3-x$, $y^3-y$, $z^3-z$, $u^2-u$, $v^2-v$ & $(1,4,8,9)$ & 0.02 & 53.00 & 108,900 \\
\hline
$x+y+z-u-v$, $x^3-x$, $y^3-y$, $z^3-z$, $u^3-u$, $v^2-v$ & $(1,4,9,12,9)$ & 0.08 & 300.00* & > 1,349,154 \\
\hline
$x+y+z-u-v+a$, $x^2-x$, $y^2-y$, $z^2-z$, $u^2-u$, $v^2-v$, $a^2-a$ & $(1,5,9) $ & < 0.01 & 8.68 & 30,030 \\
\hline
\end{tabular}
\medskip
\caption{\label{tab:results}Computational results.  The first column contains
  the considered polynomial system. The second column contains the degree
  vector of the order ideal, i.e.,
  ${(\dim \left. I^{\leq i} \middle/ I^{\leq i-1} \right.)}_i$ starting
  with $i=0$ and \(I^{\leq -1} \coloneqq 0\).
  The third column contains the average
  time (in seconds) needed to optimize a random weight over the order
  ideal polytope (we performed 20 runs for each system). The fourth column
  contains the time (in seconds) needed to count all admissible
  degree-compatible order ideals and the last column contains the actual
  number of admissible degree-compatible order ideals. The \lq{}*\rq{}
  indicates that the counting had been stopped after 300 seconds. The number
  of order ideals reported in this case is the number that have been counted
  up to that point in time.}
\end{center}
\end{table}

\section{\label{sec:Concluding-remarks}Concluding remarks}

We gave a polyhedral characterization of all order ideals that support
a border basis of a given zero-dimensional ideal.
While it is impossible
to provide a full linear description of polynomial size of the integral hull contained
in the order ideal polytope due
to Theorem~\ref{thm:xc-admissible} it might be
possible to obtain a compact extended formulation of the
order ideal polytope itself (not its integral hull). We leave this as an open question.

\bibliographystyle{plainurl}
\bibliography{gcNotesBib0,bibs}

\end{document}